# STABILIZED LOWEST ORDER FINITE ELEMENT APPROXIMATION FOR LINEAR THREE-FIELD POROELASTICITY

LORENZ BERGER,[*] RAFEL BORDAS,[†] DAVID KAY[‡] AND SIMON TAVENER[§]

**Abstract.** A stabilized conforming mixed finite element method for the three-field (displacement, fluid flux and pressure) poroelasticity problem is developed and analyzed. We use the lowest possible approximation order, namely piecewise constant approximation for the pressure and piecewise linear continuous elements for the displacements and fluid flux. By applying a local pressure jump stabilization term to the mass conservation equation we ensure stability and avoid pressure oscillations. Importantly, the discretization leads to a symmetric linear system. For the fully discretized problem we prove existence and uniqueness, an energy estimate and an optimal a-priori error estimate, including an error estimate for the divergence of the fluid flux. Numerical experiments in 2D and 3D illustrate the convergence of the method, show the effectiveness of the method to overcome spurious pressure oscillations, and evaluate the added mass effect of the stabilization term.

**1. Introduction.** Poroelasticity is a mixture theory in which a complex fluid-structure interaction is approximated by a superposition of the solid and fluid components. Developments of the continuum theory can be found, for example, in [6] and [13]. Poroelastic models have been developed to study numerous geomechanical applications ranging from reservoir engineering [30] to earthquake fault zones [36]. Fully saturated, incompressible poroelastic models have been proposed for a variety of biological tissues and processes, including lung parenchyma [23], protein-based hydrogels embedded within cells [18], blood flow in the beating myocardium [10, 12], brain oedema and hydrocephalus [25, 37], and interstitial fluid and tissue in articular cartilage and intervertebral discs [17, 20, 28].

We develop a stabilized, low-order, mixed finite element method for poroelastic models of biological tissues and restrict our attention to the fully saturated incompressible case. In order to simplify our presentation, we also assume small deformations. In contrast to [29, 36], who study a reduced displacement and pressure formulation, we retain the fluid flux variable as a primary variable resulting in a three-field, displacement, fluid flux, and pressure formulation. This avoids post processing to calculate the fluid flux and material stress, and allows physically meaningful boundary conditions to be applied at the interface when modelling the interaction between a fluid and a poroelastic structure [3]. A three-field approach can be readily extended from a Darcy to a Brinkman flow model, for which there are numerous applications in modelling biological tissues [21]. Our mixed scheme uses the lowest possible approximation order, piecewise constant approximation for the pressure and piecewise linear continuous elements for the displacement and fluid flux, since continuous pressure elements often struggle to capture the steep gradients at the interface between regions with high and low permeabilities. The resulting linear system is symmetric

[*]Department of Computer Science, University of Oxford, Wolfson Building, Parks Road, OX1 3QD, Oxford, UK, (lorenz.berger@comlab.ox.ac.uk), phone: +44 1865 273838, fax: +44 1865 273839. Funded by the EPSRC via the Life Sciences Interface Doctoral Training Centre, Oxford University.

[†]Department of Computer Science, University of Oxford, Wolfson Building, Parks Road, OX1 3QD, Oxford, UK, (rafel.bordas@comlab.ox.ac.uk), phone: +44 1865 273xxx, fax: +44 1865 273839. Funded by the EU FP7 AirPROM project (grant agreement no. 270194)

[‡]Department of Computer Science, University of Oxford, Wolfson Building, Parks Road, OX1 3QD, Oxford, UK, (david.kay@comlab.ox.ac.uk), phone: +44 1865 610814, fax: +44 1865 273839

[§]Department of Mathematics, 115 Weber Building, Colorado State University, Fort Collins, CO 80523, USA, (tavener@math.colostate.edu), phone: +1 970 491 6645, fax: +1 970 491 2161

and has a block structure that is well suited for effective preconditioning.

To ensure stability, a mixed finite element method must satisfy the Ladyzhenskaya-Babuska-Brezzi (LBB) condition. Stabilization techniques have been proposed for the Stokes equations, see e.g. [15] (chapter 5) and for Darcy flow, see e.g. [5]. Most stabilization techniques construct a modified variational formulation in which an additional term is added to the mass balance equation. In this work we use a local pressure jump stabilization method pioneered by [8] for the study of Stokes and Darcy flows that are coupled via an interface. This approach provides the natural $H^1$ stability for the displacements and $Hdiv$ stability for the fluid flux.

In earlier approaches to solving the three-field problem, [30, 31] developed a mixed finite element method using continuous piecewise linear approximations for displacements and mixed low-order Raviart Thomas elements for the fluid flux and pressure variables. However, their method was found to be susceptible to pressure oscillations [32]. To overcome these pressure oscillations, [24] and [38] analysed discontinuous and nonconforming three-field methods, respectively. Stabilization using the time derivative of pressure is shown to be crucial for stability and optimal convergence with refinement and counterexamples are provided in Section 6. In addition to these monolithic approaches there has been considerable work on operating splitting (iterative) approaches for solving the poroelastic equations [16, 22, 35]. These schemes are often only conditionally stable, and their accuracy is complicated by the exchange of information between the components, see [9].

In Section 2 we present the model and its continuous weak formulation and construct a fully-discrete approximation. We prove existence and uniqueness of solutions to this discrete model at each time step in Section 3, provide an energy estimate over time in Section 4, and derive an optimal order a-priori error estimate in Section 5. Finally in Section 6, we present numerical experiments to illustrate the convergence of the method and its ability to overcome pressure oscillations.

**2. The poroelastic model.**

**2.1. Governing equations.** Following [30] and [33], the governing equations for a fully saturated, incompressible poroelastic flow are

$$-(\lambda + \mu)\nabla(\nabla \cdot \mathbf{u}) - \mu\nabla^2\mathbf{u} + \alpha\nabla p = \mathbf{f} \quad \text{in } \Omega \times (0, T], \tag{2.1a}$$

$$k^{-1}\mathbf{z} + \nabla p = \mathbf{b} \quad \text{in } \Omega \times (0, T], \tag{2.1b}$$

$$c_0 p_t + \nabla \cdot (\alpha \mathbf{u}_t + \mathbf{z}) = g \quad \text{in } \Omega \times (0, T], \tag{2.1c}$$

$$\mathbf{u} = \mathbf{u}_D \quad \text{on } \Gamma_D \times (0, T], \tag{2.1d}$$

$$\sigma\mathbf{n} = \mathbf{t}_N \quad \text{on } \Gamma_N \times (0, T], \tag{2.1e}$$

$$\mathbf{z} \cdot \mathbf{n} = q_D \quad \text{on } \Gamma_F \times (0, T], \tag{2.1f}$$

$$p = p_D \quad \text{on } \Gamma_P \times (0, T], \tag{2.1g}$$

$$\mathbf{u}(0, \cdot) = \mathbf{u}^0 \quad \text{in } \Omega, \tag{2.1h}$$

$$p(0, \cdot) = p^0 \quad \text{in } \Omega, \tag{2.1i}$$

where $\mathbf{u}$ is the displacement, $\mathbf{z}$ is the fluid flux and $p$ is the pressure. Here $\mathbf{f}$ is the body force on the solid, $\mathbf{b}$ is the body force on the fluid and $g$ is the fluid source term, $\lambda$ and $\mu$ are the first and second Lamé parameters respectively, and $k$ is the dynamic permeability tensor. We will assume the Biot-Willis constant $\alpha = 1$, and the constrained specific storage coefficient $c_0 = 0$. Our analysis can be extended to the situation when these last two assumptions do not hold, but it is necessarily more complicated. We consider $\Omega$ to be a bounded domain in $\mathbb{R}^2$ or $\mathbb{R}^3$, and for the purpose of defining boundary conditions, $\partial\Omega = \Gamma_D \cup \Gamma_N$ for displacement and stress boundary conditions and $\partial\Omega = \Gamma_P \cup \Gamma_F$ for pressure and flux boundary conditions, with outward pointing unit normal $\mathbf{n}$.

REMARK 2.1. *Since the above resulting system of equations is linear, for ease of presentation, we will assume all Dirichlet boundary conditions are homogeneous, ie., $\mathbf{u}_D = \mathbf{0}, q_D = 0, p_D = 0$.*

**2.2. Weak formulation.** We define the following spaces for displacement, fluid flux and pressure respectively,

$$\mathbf{W}^E(\Omega) := \{\mathbf{u} \in (H^1(\Omega))^d : \mathbf{u} = \mathbf{0} \text{ on } \Gamma_D\},$$
$$\mathbf{W}^D(\Omega) := \{\mathbf{z} \in H_{div}(\Omega) : \mathbf{z} \cdot \mathbf{n} = 0 \text{ on } \Gamma_F\},$$
$$\mathcal{L}(\Omega) := \left\{ \begin{array}{ll} L^2(\Omega) & \text{if } \Gamma_N \cup \Gamma_P \neq \emptyset \\ L_0^2(\Omega) & \text{if } \Gamma_N \cup \Gamma_P = \emptyset \end{array} \right\},$$

where $L_0^2(\Omega) = \{q \in L^2(\Omega) : \int_\Omega q \, \mathrm{d}x = 0\}$. By combining these spaces we construct the mixed solution space

$$\mathcal{W}^X := \left\{\mathbf{W}^E(\Omega) \times \mathbf{W}^D(\Omega) \times \mathcal{L}(\Omega)\right\}.$$

We define the bilinear form

$$a(\mathbf{u}, \mathbf{v}) := \int_\Omega 2\mu(\epsilon(\mathbf{u}) : \epsilon(\mathbf{v})) + \lambda(\nabla \cdot \mathbf{u})(\nabla \cdot \mathbf{v}) \, \mathrm{d}x,$$

for $\mathbf{u}, \mathbf{v} \in \mathbf{W}^E(\Omega)$. This bilinear form is continuous such that

$$a(\mathbf{u}, \mathbf{v}) \leq C_c \|\mathbf{u}\|_{1,\Omega} \|\mathbf{v}\|_{1,\Omega} \quad \forall \, \mathbf{u}, \mathbf{v} \in (H^1(\Omega))^d. \tag{2.2}$$

Using Korn's inequality [7, 11], and $\int_\Omega \lambda(\nabla \cdot \mathbf{v})(\nabla \cdot \mathbf{v}) \geq 0$, we have

$$\|\mathbf{v}\|_{a,\Omega}^2 := a(\mathbf{v}, \mathbf{v}) \geq 2\mu \|\epsilon(\mathbf{v})\|_{0,\Omega}^2 \geq C_k \|\mathbf{v}_h\|_{1,\Omega}^2 \quad \forall \, \mathbf{v} \in \mathbf{W}^E(\Omega). \tag{2.3}$$

Since $k$ is assumed to be a symmetric and strictly positive definite tensor, there exists eigenfunctions $\lambda_{min}, \lambda_{max} > 0$ such that $\forall \mathbf{x} \in \Omega$, $\lambda_{min}\|\eta\|_{0,\Omega} \leq \eta^t k(\mathbf{x})\eta \leq \lambda_{max}\|\eta\|_{0,\Omega}$ $\forall \eta \in \mathbb{R}^d$, and

$$\lambda_{min}^{-1}\|\mathbf{w}\|_{0,\Omega}^2 \geq (k^{-1}\mathbf{w}, \mathbf{w}) \geq \lambda_{max}^{-1}\|\mathbf{w}\|_{0,\Omega}^2 \quad \forall \mathbf{w} \in \mathbf{W}^D(\Omega). \tag{2.4}$$

The continuous weak problem is: Find $\mathbf{u} \in \mathbf{W}^E(\Omega)$, $\mathbf{z} \in \mathbf{W}^D(\Omega)$, and $p \in \mathcal{L}(\Omega)$ for any time $t \in [0, T]$ such that

$$a(\mathbf{u}, \mathbf{v}) - (p, \nabla \cdot \mathbf{v}) = (\mathbf{f}, \mathbf{v}) + (\mathbf{t}_N, \mathbf{v})_{\Gamma_N} \quad \forall \mathbf{v} \in \mathbf{W}^E(\Omega), \tag{2.5a}$$

$$(k^{-1}\mathbf{z}, \mathbf{w}) - (p, \nabla \cdot \mathbf{w}) = (\mathbf{b}, \mathbf{w}) \quad \forall \mathbf{w} \in \mathbf{W}^D(\Omega), \tag{2.5b}$$

$$(\nabla \cdot \mathbf{u}_t, q) + (\nabla \cdot \mathbf{z}, q) = (g, q) \quad \forall q \in \mathcal{L}(\Omega). \tag{2.5c}$$

We will assume the following regularity requirements on the data,

$$\begin{aligned} \mathbf{f} &\in C^1([0,T]; (H^{-1}(\Omega))^d), & \mathbf{t}_N &\in C^1([0,T]; H^{-1/2}(\Gamma_N)), \\ \mathbf{b} &\in C^1([0,T]; H_{div}^{-1}(\Omega)), & g &\in C^0([0,T]; (L^2(\Omega))^d). \end{aligned} \tag{2.6}$$

For the initial conditions we require that $\mathbf{u}^0 \in (H^1(\Omega))^d$. The well-posedness of the continuous two-field formulation has been proven by [33]. [26] proves well-posedness for the continuous three-field formulation (2.5). In this work we also establish the well-posedness of (2.5) as a result of the energy estimates proven in section 4, see remark 4.1.

**2.3. Fully-discrete model.** Let $\mathcal{T}^h$ be a partition of $\Omega$ into non-overlapping elements $K$, where $h$ denotes the size of the largest element in $\mathcal{T}^h$ and assume that the partition is quasi-uniform. We define the following finite element spaces,

$$\mathbf{W}_h^E := \left\{ \mathbf{u}_h \in C^0(\Omega) : \mathbf{u}_h|_K \in P_1(K) \ \forall K \in \mathcal{T}^h, \mathbf{u}_h = 0 \text{ on } \Gamma_D \right\},$$
$$\mathbf{W}_h^D := \left\{ \mathbf{z}_h \in C^0(\Omega) : \mathbf{z}_h|_K \in P_1(K) \ \forall K \in \mathcal{T}^h, \mathbf{z}_h \cdot \mathbf{n} = 0 \text{ on } \Gamma_F \right\},$$
$$Q_h := \begin{cases} \{p_h : p_h|_K \in P_0(K) \ \forall K \in \mathcal{T}^h\} & \text{if } \Gamma_N \cup \Gamma_P \neq \emptyset \\ \{p_h : p_h|_K \in P_0(K), \int_\Omega p_h = 0 \ \forall K \in \mathcal{T}^h\} & \text{if } \Gamma_N \cup \Gamma_P = \emptyset \end{cases},$$

where $P_0(K)$ and $P_1(K)$ are respectively the spaces of constant and linear polynomials on $K$. We partition $[0,T]$ into $N$ evenly spaced non-overlapping regions $(t_{n-1}, t_n]$, $n = 1, 2, \ldots, N$, where $t_n - t_{n-1} = \Delta t$. For any sufficiently smooth function $v(t,x)$ we define $v^n(x)" = v(t_n, x)$ and the discrete time derivative by $v^n_{\Delta t} := \frac{v^n - v^{n-1}}{\Delta t}$.

The fully discrete weak problem is: For $n = 1, 2, \ldots, N$, $\mathbf{u}_h^n \in \mathbf{W}_h^E$, find $\mathbf{z}_h^n \in \mathbf{W}_h^D$ and $p_h^n \in Q_h$ such that

$$a(\mathbf{u}_h^n, \mathbf{v}_h) - (p_h^n, \nabla \cdot \mathbf{v}_h) = (\mathbf{f}^n, \mathbf{v}_h) + (\mathbf{t}_N, \mathbf{v}_h)_{\Gamma_N} \ \forall \mathbf{v}_h \in \mathbf{W}_h^E, \quad (2.7a)$$

$$(k^{-1} \mathbf{z}_h^n, \mathbf{w}_h) - (p_h^n, \nabla \cdot \mathbf{w}_h) = (\mathbf{b}^n, \mathbf{w}_h) \ \forall \mathbf{w}_h \in \mathbf{W}_h^D, \quad (2.7b)$$

$$(\nabla \cdot \mathbf{u}_{\Delta t, h}, q_h) + (\nabla \cdot \mathbf{z}_h^n, q_h) + J(p_{\Delta t, h}, q_h) = (g^n, q_h) \ \forall q_h \in Q_h. \quad (2.7c)$$

The stabilization term is

$$J(p,q) = \delta \sum_K \int_{\partial K \setminus \partial \Omega} h_{\partial K} [p][q] \, ds, \quad (2.8)$$

where $h_{\partial K}$ denotes the size (diameter) of an element edge in 2D or face in 3D, $\delta$ is a stabilization parameter that is independent of $h$ and $\Delta t$, and $[\cdot]$ is the jump across an edge or face (taken on the interior edges only).

We also assume

$$a(\mathbf{u}_h^0, \mathbf{v}_h) = a(\mathbf{u}^0, \mathbf{v}_h) \ \forall \mathbf{v}_h \in \mathbf{W}_h^E, \quad (2.9a)$$

$$J(p_h^0, q_h) = J(p^0, q_h) \ \forall q_h \in Q_h, \quad (2.9b)$$

where $p^0 \in \mathcal{L}(\Omega)$.

**3. Existence and uniqueness of solutions to the fully-discrete model.** Well-posedness of the unstabilized fully-discretized system (2.7), $\delta = 0$, with the use of a low order Raviart-Thomas approximation for the fluid velocity is shown by [31] for $c_0 > 0$, and by [26] for $c_0 \geq 0$. Although as the permeability tends to zero and the porous mixture becomes impermeable, the three-field linear poroelasticity tends to a mixed linear elasticity problem [19]. Hence, in this case this element becomes unstable, as expected since the elasticity $P1 - P0$ approximation is known to be unstable. Our method is stable for both the Darcy problem (as the elasticity coefficients tend to infinity) and the mixed linear elasticity problem (as the permeability tends to zero), and is therefore stable for all permeabilities and elasticity coefficients.

**3.1. Norms and inequalities.** The stabilization term gives rise to the semi-norm $|q|_{J,\Omega} := J(q,q)^{1/2}$. Throughout this work, we will let $C$ denote a generic positive constant, whose value may change from instance to instance, but is independent of any mesh parameters. Using the scaling argument $\left\| h^{1/2} p_h \right\|_{0, \partial K} \leq c_z \|p_h\|_{0,K}$,

Cauchy-Schwarz and the triangle inequality the following bounds for the stabilization term hold.

$$|p_h|_{J,\Omega} \leq C\|p_h\|_{0,\Omega} \text{ and } J(p_h, q_h) \leq |p_h|_{J,\Omega}|q_h|_{J,\Omega} \; \forall \; p_h, q_h \in Q_h. \tag{3.1}$$

Furthermore, for any $q \in H^1(\Omega)$,

$$J(p, q) = 0 \; \forall p \in \mathcal{L}(\Omega), \tag{3.2}$$

see Lemma 1.23 in [14]. We also have the Poincaré inequality

$$\|q\|_{0,\Omega} \leq C_p \|\nabla q\|_{0,\Omega} \; \; \forall q \in H^1(\Omega).$$

We now give some approximation results that will be useful later. Let $\pi_h^1 : H^1(\Omega) \to \mathbf{W}_h^E$ and $\pi_h^0 : L^2(\Omega) \to Q_h$ be Clément projections, see [11].

LEMMA 3.1. *For all* $\mathbf{v} \in \left(H^2(\Omega)\right)^d$ *and* $q \in H^1(\Omega)$ *the interpolation operators satisfy: For* $s = 0, 1$

$$\|\mathbf{v} - \pi_h^1 \mathbf{v}\|_{s,\Omega} \leq Ch^{2-s}\|\mathbf{v}\|_{2,\Omega}, \tag{3.3}$$

$$\|q - \pi_h^0 q\|_{0,\Omega} \leq Ch\|q\|_{1,\Omega}, \tag{3.4}$$

$$|q - \pi_h^0 q|_{J,\Omega} \leq Ch\|q\|_{1,\Omega}. \tag{3.5}$$

*Proof.* The first two results are standard. The final result is obtained by using the element error estimate provided in [34] and then summing over all elements. □

Due to the surjectivity of the divergence operator, for every $p \in L^2(\Omega)$ there exists a function $\mathbf{v}_p \in (H_0^1(\Omega))^d$ such that $\nabla \cdot \mathbf{v}_p = -p$ and $\|\mathbf{v}_p\|_{1,\Omega} \leq c\|p\|_{0,\Omega}$. We assume that the projection, $\pi_h^1 \mathbf{v}_p \in \mathbf{v}_p \in (H_0^1(\Omega))^d$, is stable such that

$$\|\pi_h^1 \mathbf{v}_p\|_{1,\Omega} \leq \hat{c}\|p\|_{0,\Omega}. \tag{3.6}$$

Furthermore, for any element $K \in \mathcal{T}^h$

$$\|\mathbf{v}_p - \pi_h^1 \mathbf{v}_p\|_{0,K} \leq Ch\|\mathbf{v}_p\|_{H^1(\omega_K)}, \tag{3.7}$$

where $\omega_K$ is the union of all elements $J \in \mathcal{T}^h$ such that $\overline{K} \cap \overline{J} \neq \emptyset$.

Combining the above with the trace inequality, see Lemma 3.1 in [34],

$$\left\|(\mathbf{v}_p - \pi_h^1 \mathbf{v}_p) \cdot \mathbf{n}\right\|_{0,\partial K}^2 \leq C\|\mathbf{v}_p - \pi_h^1 \mathbf{v}_p\|_{0,K}(h^{-1}\|\mathbf{v}_p - \pi_h^1 \mathbf{v}_p\|_{0,K} + \|\mathbf{v}_p - \pi_h^1 \mathbf{v}_p\|_{1,K}), \tag{3.8}$$

we obtain

$$\left\|(\mathbf{v}_p - \pi_h^1 \mathbf{v}_p) \cdot \mathbf{n}\right\|_{0,\partial K}^2 \leq Ch\|\mathbf{v}_p\|_{H^1(\omega_K)}^2. \tag{3.9}$$

Taking into account $\|\mathbf{v}_p\|_{1,\Omega} \leq c\|p\|_{0,\Omega}$, we may write

$$\sum_K \int_{\partial K} h^{-1}|(\mathbf{v}_p - \pi_h^1 \mathbf{v}_p) \cdot \mathbf{n}|^2 \, \mathrm{d}s \leq c_t\|p\|_{0,\Omega}^2. \tag{3.10}$$

We define the fully discrete finite element approximation for all time to be the piecewise constant in time functions $\mathbf{u}_h(t, \mathbf{x}) := \mathbf{u}_h^n(\mathbf{x})$ for $t \in (t_{n-1}, t_n]$, $\mathbf{z}_h$ and $p_h$ are defined similarly. For such piecewise continuous in time functions, $v$, the norms $L^2(0,T;X)$ satisfy

$$\|v\|_{(L^2;X)}^2 = \int_0^T \|v(s,\cdot)\|_X^2 \, ds = \sum_{n=1}^N \Delta t \|v^n\|_X^2,$$

where $X$ is any given function space over $\Omega$.

For all $\mathbf{v} \in H^2(0,T;(L^2(\Omega))^d)$

$$\sum_{n=1}^N \Delta t \left\| \mathbf{v}_{\Delta t}^n - \frac{\partial \mathbf{v}}{\partial t}(t^n,\cdot) \right\|_{0,\Omega}^2 \leq \Delta t^2 \int_0^T \|\mathbf{v}_{tt}\|_{0,\Omega}^2 ds. \qquad (3.11)$$

For all $[\mathbf{v},\mathbf{w},q] \in \left[(H^1(\Omega))^d \times H_{div}(\Omega) \times L^2(\Omega)\right]$ we define the norm

$$\||[\mathbf{v},\mathbf{w},q]\||_A^2 := \|\mathbf{v}\|_{1,\Omega}^2 + \Delta t^2 \|\nabla \cdot \mathbf{w}\|_{0,\Omega}^2 + \Delta t \|\mathbf{w}\|_{0,\Omega}^2 + \|q\|_{0,\Omega}^2 + |q|_{J,\Omega}^2, \qquad (3.12)$$

and for all $[\mathbf{v},\mathbf{w},q] \in \left[L^\infty(0,T;(H^1(\Omega))^d) \times L^2(0,T;H_{div}(\Omega)) \times L^2(0,T;L^2(\Omega))\right]$ the norm

$$\||[\mathbf{v},\mathbf{w},q]\||_B^2 := \|\mathbf{v}\|_{L^\infty(H^1)}^2 + \|\mathbf{w}\|_{L^2(L^2)}^2 + \|q\|_{L^2(L^2)}^2. \qquad (3.13)$$

**3.2. Existence and uniqueness.** Combining the fully discrete equations (2.7a), (2.7b) and (2.7c), after first multiplying (2.7b) and (2.7c) by $\Delta t$, gives the equivalent problem: For $n = 1, 2, \ldots, n$ find $(\mathbf{u}_h, \mathbf{z}_h, p_h)$ such that

$$\begin{aligned} B_h^n[(\mathbf{u}_h,\mathbf{z}_h,p_h),(\mathbf{v}_h,\mathbf{w}_h,q_h)] \\ = (\mathbf{f}^n,\mathbf{v}_h) + (\mathbf{t}_N,\mathbf{v}_h)_{\Gamma_N} + \Delta t(\mathbf{b}^n,\mathbf{w}_h) + \Delta t(g^n,q_h) \\ + (\nabla \cdot \mathbf{u}_h^{n-1},q_h) + J(p_h^{n-1},q_h) \quad \forall (\mathbf{v}_h,\mathbf{w}_h,q_h) \in \mathcal{W}_h^X, \end{aligned}$$

where

$$\begin{aligned} B_h^n[(\mathbf{u}_h,\mathbf{z}_h,p_h),(\mathbf{v}_h,\mathbf{w}_h,q_h)] \\ = a(\mathbf{u}_h^n,\mathbf{v}_h) + \Delta t(k^{-1}\mathbf{z}_h^n,\mathbf{w}_h) - (p_h^n,\nabla \cdot \mathbf{v}_h) - \Delta t(p_h^n,\nabla \cdot \mathbf{w}_h) \\ + (\nabla \cdot \mathbf{u}_h^n,q_h) + \Delta t(\nabla \cdot \mathbf{z}_h^n,q_h) + J(p_h^n,q_h). \quad (3.14) \end{aligned}$$

The linear form satisfies the following continuity property

$$|B_h^n[(\mathbf{u}_h,\mathbf{z}_h,p_h),(\mathbf{v}_h,\mathbf{w}_h,q_h)]| \leq C \,\||(\mathbf{u}_h^n,\mathbf{z}_h^n,p_h^n)\||_A \, \||(\mathbf{v}_h,\mathbf{w}_h,q_h)\||_A \,.$$

We apply Babuska's theory [2] to show well-posedness (existence and uniqueness) of this discretized system at a particular time step. This requires us to prove a discrete inf-sup type result (Theorem 3.2) for the combined bilinear form (3.14).

THEOREM 3.2. *Let $\gamma > 0$ be a constant independent of any mesh parameters. Then the finite element formulation (2.7) satisfies the following discrete inf-sup condition*

$$\gamma \,\||(\mathbf{u}_h^n,\mathbf{z}_h^n,p_h^n)\||_A \leq \sup_{(\mathbf{v}_h,\mathbf{w}_h,q_h) \in \mathcal{V}_h^X} \frac{B_h^n[(\mathbf{u}_h,\mathbf{z}_h,p_h),(\mathbf{v}_h,\mathbf{w}_h,q_h)]}{\||(\mathbf{v}_h,\mathbf{w}_h,q_h)\||_A} \quad \forall (\mathbf{u}_h,\mathbf{z}_h,p_h) \in \mathcal{W}_h^X.$$

$$(3.15)$$

Hence, given a solution at the previous time step the linear system arising from the fully discrete method for the subsequent time step is non-singular. The following proof follows ideas presented by [8].

*Proof.*

*Step 1, bounding* $\|\mathbf{u}_h^n\|_{1,\Omega}$, $\Delta t^{1/2}\|\mathbf{z}_h^n\|_{0,\Omega}$, *and* $|p_h^n|_{J,\Omega}$.
Choose $(\mathbf{v}_h, \mathbf{w}_h, q_h) = (\mathbf{u}_h^n, \mathbf{z}_h^n, p_h^n)$, then using (2.3) and (2.4), we obtain

$$B_h^n[(\mathbf{u}_h, \mathbf{z}_h, p_h), (\mathbf{u}_h, \mathbf{z}_h, p_h)] = a(\mathbf{u}_h^n, \mathbf{u}_h^n) + \Delta t(k^{-1}\mathbf{z}_h^n, \mathbf{z}_h^n) + J(p_h^n, p_h^n)$$
$$\geq C_k \|\mathbf{u}_h^n\|_{1,\Omega}^2 + \lambda_{max}^{-1}\Delta t\|\mathbf{z}_h^n\|_{0,\Omega}^2 + |p_h^n|_{J,\Omega}^2. \quad (3.16)$$

*Step 2, bounding* $\|p_h^n\|_{0,\Omega}$.
Choose $(\mathbf{v}_h, \mathbf{w}_h, q_h) = (\pi_h^1 \mathbf{v}_{p_h^n}, 0, 0)$ and add $0 = \|p_h^n\|_{0,\Omega}^2 + (p_h^n, \nabla \cdot \mathbf{v}_{p_h^n})$ to obtain

$$B_h^n[(\mathbf{u}_h, \mathbf{z}_h, p_h), (\pi_h^1 \mathbf{v}_{p_h^n}, 0, 0)] = a(\mathbf{u}_h^n, \pi_h^1 \mathbf{v}_{p_h^n})$$
$$+ \|p_h^n\|_{0,\Omega}^2 + (p_h^n, \nabla \cdot (\mathbf{v}_{p_h^n} - \pi_h^1 \mathbf{v}_{p_h^n})). \quad (3.17)$$

Focusing on the third term in (3.17) only, by applying the divergence theorem together with the fact that $p_h^n$ is piecewise constant on $\mathcal{T}^h$, and recalling $\mathbf{v}_{p_h^n} - \pi_h^1 \mathbf{v}_{p_h^n} \in (H_0^1(\Omega))^d$, we obtain

$$(p_h^n, \nabla \cdot (\mathbf{v}_{p_h^n} - \pi_h^1 \mathbf{v}_{p_h^n})) = \sum_K \int_{\partial K/\partial\Omega} p_h^n (\mathbf{v}_{p_h^n} - \pi_h^1 \mathbf{v}_{p_h^n}) \cdot \mathbf{n}\, ds$$
$$= \sum_K \frac{1}{2} \int_{\partial K/\partial\Omega} [p_h^n](\mathbf{v}_{p_h^n} - \pi_h^1 \mathbf{v}_{p_h^n}) \cdot \mathbf{n}\, ds.$$

We thus have

$$B_h^n[(\mathbf{u}_h, \mathbf{z}_h, p_h), (\pi_h^1 \mathbf{v}_{p_h^n}, 0, 0)] = \|p_h^n\|_{0,\Omega}^2 + a(\mathbf{u}_h^n, \pi_h^1 \mathbf{v}_{p_h^n})$$
$$+ \sum_K \frac{1}{2} \int_{\partial K/\partial\Omega} [p_h^n](\mathbf{v}_{p_h^n} - \pi_h^1 \mathbf{v}_{p_h^n}) \cdot \mathbf{n}\, ds.$$

Now first applying the Cauchy-Schwarz inequality and (2.2) on the right hand side to get

$$B_h^n[(\mathbf{u}_h, \mathbf{z}_h, p_h), (\pi_h^1 \mathbf{v}_{p_h^n}, 0, 0)] \geq \|p_h^n\|_{0,\Omega}^2 - C_c \|\mathbf{u}_h^n\|_{1,\Omega} \|\pi_h^1 \mathbf{v}_{p_h^n}\|_{1,\Omega}$$
$$- \sum_K \frac{1}{2} \left( \int_{\partial K/\partial\Omega} \left(h^{1/2}[p_h^n]\right)^2 ds \right)^{1/2} \cdot \left( \int_{\partial K} \left(h^{-1/2}(\mathbf{v}_{p_h^n} - \pi_h^1 \mathbf{v}_{p_h^n}) \cdot \mathbf{n}\right)^2 ds \right)^{1/2}.$$

Now apply Young's inequality and (3.6) to obtain

$$B_h^n[(\mathbf{u}_h, \mathbf{z}_h, p_h), (\pi_h^1 \mathbf{v}_{p_h^n}, 0, 0)] \geq \|p_h^n\|_{0,\Omega}^2 - \frac{C_c^2}{2\epsilon} \|\mathbf{u}_h^n\|_{1,\Omega}^2 - \frac{\epsilon \hat{c}}{2} \|p_h^n\|_{0,\Omega}^2$$
$$- \frac{1}{2\epsilon\delta} J(p_h^n, p_h^n) - \frac{\epsilon}{2} \sum_K \int_{\partial K} h^{-1} |(\mathbf{v}_{p_h^n} - \pi_h^1 \mathbf{v}_{p_h^n}) \cdot \mathbf{n}|^2 ds.$$

Applying (3.10) we obtain

$$B_h^n[(\mathbf{u}_h, \mathbf{z}_h, p_h), (\pi_h^1 \mathbf{v}_{p_h^n}, 0, 0)] \geq -\frac{C_c^2}{2\epsilon}\|\mathbf{u}_h^n\|_{1,\Omega}^2 + \left(1 - (\hat{c} + c_t)\frac{\epsilon}{2}\right)\|p_h^n\|_{0,\Omega}^2 - \frac{1}{2\epsilon\delta}|p_h^n|_{J,\Omega}^2. \quad (3.18)$$

*Step 3, bounding* $\Delta t \|\nabla \cdot \mathbf{z}_h^n\|_{0,\Omega}$.
Choosing $(\mathbf{v}_h, \mathbf{w}_h, q_h) = (0, 0, \Delta t \nabla \cdot \mathbf{z}_h^n)$ yields

$$B_h^n[(\mathbf{u}_h, \mathbf{z}_h, p_h), (0, 0, \Delta t\nabla \cdot \mathbf{z}_h^n)] = (\nabla \cdot \mathbf{u}_h^n, \Delta t\nabla \cdot \mathbf{z}_h^n) + \Delta t^2 \|\nabla \cdot \mathbf{z}_h^n\|_{0,\Omega}^2 + J(p_h^n, \Delta t\nabla \cdot \mathbf{z}_h^n).$$

We bound the first term using the Cauchy-Schwarz inequality followed by Young's inequality such that

$$(\nabla \cdot \mathbf{u}_h^n, \Delta t\nabla \cdot \mathbf{z}_h^n) \leq \frac{C_p}{2\epsilon}\|\mathbf{u}_h^n\|_{1,\Omega}^2 + \frac{\epsilon\Delta t^2}{2}\|\nabla \cdot \mathbf{z}_h^n\|_{0,\Omega}^2.$$

We can also bound the third term as before using the Cauchy-Schwarz inequality followed by Young's inequality such that

$$\begin{aligned} J(p_h^n, \Delta t\nabla \cdot \mathbf{z}_h^n) &\leq \frac{1}{2\epsilon}J(p_h^n, p_h^n) + \frac{\epsilon\Delta t^2}{2}J(\nabla \cdot \mathbf{z}_h^n, \nabla \cdot \mathbf{z}_h^n) \\ &= \frac{1}{2\epsilon}J(p_h^n, p_h^n) + \epsilon\delta\Delta t^2 \sum_K \int_{\partial K/\partial\Omega} |h^{1/2}\nabla \cdot \mathbf{z}_h^n|^2 \, \mathrm{d}s \\ &\leq \frac{1}{2\epsilon}J(p_h^n, p_h^n) + \epsilon\delta c_z \Delta t^2 \|\nabla \cdot \mathbf{z}_h^n\|_{0,\Omega}^2, \end{aligned} \quad (3.19)$$

where we have used the fact that $\nabla \cdot \mathbf{z}_h^n \in Q_h$ in conjunction with (3.1). This yields

$$B_h^n[(\mathbf{u}_h, \mathbf{z}_h, p_h), (0, 0, \Delta t\nabla \cdot \mathbf{z}_h^n)] \geq (1 - \epsilon\delta c_z - \frac{\epsilon}{2})\Delta t^2 \|\nabla \cdot \mathbf{z}_h^n\|_{0,\Omega}^2 - \frac{1}{2\epsilon}|p_h^n|_{J,\Omega}^2 - \frac{C_p}{2\epsilon}\|\mathbf{u}_h^n\|_{1,\Omega}^2. \quad (3.20)$$

*Step 4, Combining steps 1-3.* Finally, we can combine (3.16), (3.18) and (3.20) to get control over all the norms by choosing

$$(\mathbf{v}_h, \mathbf{w}_h, q_h) = (\beta\mathbf{u}_h^n + \pi_h^1 \mathbf{v}_{p_h^n}, \beta\mathbf{z}_h^n, \beta p_h^n + \Delta t\nabla \cdot \mathbf{z}_h^n),$$

where $\beta$ is a real number that will be chosen to be sufficiently large enough to conclude, which yields

$$B_h^n[(\mathbf{u}_h, \mathbf{z}_h, p_h), (\beta\mathbf{u}_h^n + \pi_h^1\mathbf{v}_{p_h^n}, \beta\mathbf{z}_h^n, \beta p_h^n + \Delta t\nabla \cdot \mathbf{z}_h^n)] \geq$$
$$(\beta C_k - \frac{C_c^2 + C_p}{2\epsilon})\|\mathbf{u}_h^n\|_{1,\Omega}^2 + \beta\lambda_{max}^{-1}\Delta t\|\mathbf{z}_h^n\|_{0,\Omega}^2 + \left(1 - \epsilon\delta c_z - \frac{\epsilon}{2}\right)\Delta t^2\|\nabla \cdot \mathbf{z}_h^n\|_{0,\Omega}^2$$
$$+ \left(1 - (\hat{c} + c_t)\frac{\epsilon}{2}\right)\|p_h^n\|_{0,\Omega}^2 + \left(\beta - \frac{1}{2\epsilon} - \frac{1}{2\epsilon\delta}\right)|p_h^n|_{J,\Omega}^2, \quad (3.21)$$

where we can choose

$$\beta \geq \max\left[\frac{C_c^2 + C_p}{2\epsilon C_k} + \frac{1 - \bar{C}\epsilon}{C_k}, \lambda_{max}\left(1 - \bar{C}\epsilon\right), \frac{1}{2\epsilon} + \frac{1}{2\epsilon\delta} + 1 - \bar{C}\epsilon\right],$$

with $\bar{C} = \max\left[\frac{\hat{c}+c_t}{2}, \delta c_z - \frac{1}{2}\right]$. This yields

$$B_h^n[(\mathbf{u}_h, \mathbf{z}_h, p_h), (\beta\mathbf{u}_h^n + \pi_h^1 \mathbf{v}_{p_h^n}, \beta\mathbf{z}_h^n, \beta p_h^n + \nabla \cdot \mathbf{z}_h^n)] \geq (1 - \bar{C}\epsilon) |||(\mathbf{u}_h^n, \mathbf{z}_h^n, p_h^n)|||_A^2.$$

To complete the proof, we let $(\mathbf{v}_h, \mathbf{w}_h, q_h) = (\beta\mathbf{u}_h^n + \pi_h^1 \mathbf{v}_{p_h^n}, \beta\mathbf{z}_h^n, \beta p_h^n + \Delta t \nabla \cdot \mathbf{z}_h^n)$ and show that for $\epsilon$ sufficiently small there exists a constant $C$ such that $|||(\mathbf{u}_h^n, \mathbf{z}_h^n, p_h^n)|||_A \geq C |||(\mathbf{v}_h, \mathbf{w}_h, q_h)|||_A$. Using the triangle inequality and (3.6) we obtain

$$\begin{aligned}
&|||(\beta\mathbf{u}_h^n + \pi_h^1 \mathbf{v}_{p_h^n}, \beta\mathbf{z}_h^n, \beta p_h^n + \Delta t \nabla \cdot \mathbf{z}_h^n)|||_A^2 \\
&\leq C\left(\beta^2 \|\mathbf{u}_h^n\|_{1,\Omega}^2 + \|\pi_h^1 \mathbf{v}_{p_h^n}\|_{1,\Omega}^2 + \Delta t^2 (1+\beta)^2 \|\nabla \cdot \mathbf{z}_h^n\|_{0,\Omega}^2 + \beta^2 \Delta t \|\mathbf{z}_h^n\|_{0,\Omega}^2 \right.\\
&\qquad\left. + \beta^2 \|p_h^n\|_{0,\Omega}^2 + \beta^2 |p_h^n|_{J,\Omega}^2 + \Delta t^2 |\nabla \cdot \mathbf{z}_h^n|_{J,\Omega}^2\right) \\
&\leq C |||(\mathbf{u}_h^n, \mathbf{z}_h^n, p_h^n)|||_A^2,
\end{aligned}$$

as desired. □

**4. Energy estimate for the fully-discrete model.** In this Section we construct two new combined bilinear forms, $B_{\Delta t,h}^n$ (Lemmas 4.1 and 4.2) and $\mathcal{B}_h^n$ (Lemmas 4.3 and 4.4). These bilinear forms are bounded below by Lemmas 4.1 and 4.3 respectively. Lemma 4.2 uses Lemma 4.1 to provide a bound on $\mathbf{u}_h, \mathbf{z}_h$ and $p_h$. Lemma 4.4 uses Lemma 4.3 to provide a bound on $\nabla \cdot \mathbf{z}_h$.

**4.1. Bound on the displacement, fluid flux and pressure.** Adding (2.7a), (2.7b) and (2.7c), and assuming $\mathbf{t}_N = 0$ on $\Gamma_N$, we get the following

$$B_{\Delta t,h}^n[(\mathbf{u}_h, \mathbf{z}_h, p_h), (\mathbf{v}_h, \mathbf{w}_h, q_h)] = (\mathbf{f}^n, \mathbf{v}_h) + (\mathbf{b}^n, \mathbf{w}_h) + (g^n, q_h) \quad \forall (\mathbf{v}_h, \mathbf{w}_h, q_h) \in \mathcal{W}_h^X, \tag{4.1}$$

where

$$\begin{aligned}
B_{\Delta t,h}^n[(\mathbf{u}_h, \mathbf{z}_h, p_h), (\mathbf{v}_h, \mathbf{w}_h, q_h)] &= a(\mathbf{u}_h^n, \mathbf{v}_h) + (k^{-1}\mathbf{z}_h^n, \mathbf{w}_h) - (p_h^n, \nabla \cdot \mathbf{v}_h) - (p_h^n, \nabla \cdot \mathbf{w}_h) \\
&\quad + (\nabla \cdot \mathbf{u}_{\Delta t,h}^n, q_h) + (\nabla \cdot \mathbf{z}_h^n, q_h) + J(p_{\Delta t,h}^n, q_h).
\end{aligned} \tag{4.2}$$

LEMMA 4.1. $(\mathbf{u}_h, \mathbf{z}_h, p_h)$ satisfies

$$\begin{aligned}
\sum_{n=1}^N \Delta t B_{\Delta t,h}^n &[(\mathbf{u}_h, \mathbf{z}_h, p_h), (\mathbf{u}_{\Delta t,h}^n + \pi_h^1 \mathbf{v}_{p_h^n}, \mathbf{z}_h^n, p_h^n)] \\
&+ \|\mathbf{u}_h^0\|_{1,\Omega}^2 + |p_h^0|_{J,\Omega}^2 + \|\mathbf{u}_h\|_{L^2(H^1)}^2 + \|p_h\|_{L^2(J)}^2 \\
&\geq C\left(\|\mathbf{u}_h^N\|_{1,\Omega}^2 + |p_h^N|_{J,\Omega}^2 + \|\mathbf{z}_h\|_{L^2(L^2)}^2 + \|p_h\|_{L^2(L^2)}^2\right).
\end{aligned}$$

*Proof.*
For $n = 1, 2, \ldots, N$ we choose $(\mathbf{v}_h, \mathbf{w}_h, q_h) = (\mathbf{u}_{\Delta t,h}^n + \pi_h^1 \mathbf{v}_{p_h^n}, \mathbf{z}_h^n, p_h^n)$ in (4.2),

multiplying by $\Delta t$, and summing over all time steps, we get

$$\sum_{n=1}^{N} \Delta t B_{\Delta t,h}^n[(\mathbf{u}_h, \mathbf{z}_h, p_h),(\mathbf{u}_{\Delta t,h}^n + \pi_h^1 \mathbf{v}_{p_h^n}, \mathbf{z}_h^n, p_h^n)]$$
$$= \sum_{n=1}^{N} \Delta t a(\mathbf{u}_h^n, \mathbf{u}_{\Delta t,h}^n) + \sum_{n=1}^{N} \Delta t J(p_{\Delta t,h}^n, p_h^n) + \sum_{n=1}^{N} \Delta t (k^{-1} \mathbf{z}_h^n, \mathbf{z}_h^n)$$
$$+ \sum_{n=1}^{N} \Delta t a(\mathbf{u}_h^n, \pi_h^1 \mathbf{v}_{p_h^n}) - \sum_{n=1}^{N} \Delta t (p_h^n, \nabla \cdot \pi_h^1 \mathbf{v}_{p_h^n}). \quad (4.3)$$

By telescoping out the first two terms on the righthand side, using (2.4) on the third, and applying firstly Young's inequality and then (3.6) to the final two terms, we obtain the inequality

$$\left( \sum_{n=0}^{N} \Delta t B_{\Delta t,h}^n[(\mathbf{u}_h, \mathbf{z}_h, p_h),(\mathbf{u}_{\Delta t,h}^n + \pi_h^1 \mathbf{v}_{p_h^n}, \mathbf{z}_h^n, p_h^n)] \right.$$
$$\left. + \frac{C_c}{2} \|\mathbf{u}_h^0\|_{1,\Omega}^2 + \frac{C_c}{2\epsilon} \|\mathbf{u}_h\|_{L^2(H^1)}^2 + \frac{1}{4\epsilon} \|p_h\|_{L^2(J)}^2 + \frac{1}{2} |p_h^0|_{J,\Omega}^2 \right)$$
$$\geq \frac{C_k}{2} \|\mathbf{u}_h^N\|_{1,\Omega}^2 + \frac{1}{2} |p_h^N|_{J,\Omega}^2 + \lambda_{max}^{-1} \|\mathbf{z}_h\|_{L^2(L^2)}^2 + (1 - C\epsilon) \|p_h\|_{L^2(L^2)}^2. \quad (4.4)$$

Finally, choosing $\epsilon$ sufficiently small completes the proof. □

LEMMA 4.2. $(\mathbf{u}_h, \mathbf{z}_h, p_h)$ satisfies

$$\left\| \mathbf{u}_h^N \right\|_{1,\Omega}^2 + |p_h^N|_{J,\Omega}^2 + \|\mathbf{z}_h\|_{L^2(L^2)}^2 + \|p_h\|_{L^2(L^2)}^2 \leq C(T).$$

*Proof.*
For $n = 1, 2, \ldots, N$ we choose $(\mathbf{v}_h, \mathbf{w}_h, q_h) = (\mathbf{u}_{\Delta t,h}^n + \pi_h^1 \mathbf{v}_{p_h^n}, \mathbf{z}_h^n, p_h^n)$ in (4.1), multiplying by $\Delta t$, and summing yields

$$\sum_{n=1}^{N} \Delta t B_{\Delta t,h}^n[(\mathbf{u}_h^n, \mathbf{z}_h^n, p_h^n),(\mathbf{u}_{\Delta t,h}^n + \pi_h^1 \mathbf{v}_{p_h^n}, \mathbf{z}_h^n, p_h^n)] = \sum_{n=1}^{N} \Delta t (\mathbf{f}^n, \mathbf{u}_{\Delta t,h}^n + \pi_h^1 \mathbf{v}_{p_h^n})$$
$$+ \sum_{n=1}^{N} \Delta t (\mathbf{b}^n, \mathbf{z}_h^n) + \sum_{n=1}^{N} \Delta t (g^n, p_h^n).$$

Let us note the standard result, for any $\epsilon > 0$

$$\sum_{n=1}^{N} \Delta t (\mathbf{f}^n, \mathbf{u}_{\Delta t,h}^n) \leq C \left[ \frac{1}{2\epsilon} \left( \|\mathbf{f}^0\|_{0,\Omega}^2 + \|\mathbf{f}^N\|_{0,\Omega}^2 + \|\mathbf{f}_t\|_{L^2(L^2)}^2 \right) \right.$$
$$\left. + \frac{\epsilon}{2} \left( \|\mathbf{u}_h^0\|_{0,\Omega}^2 + \|\mathbf{u}_h^N\|_{0,\Omega}^2 + \|\mathbf{u}_h\|_{L^2(L^2)}^2 \right) \right]. \quad (4.5)$$

Now using the above, Lemma 4.1, the Cauchy-Schwarz and Young's inequalities,

choosing $\epsilon$ sufficiently small, and noting (3.6), we arrive at

$$\left\|\mathbf{u}_h^N\right\|_{1,\Omega}^2 + |p_h^N|_{J,\Omega}^2 + \|\mathbf{z}_h\|_{L^2(L^2)}^2 + \|p_h\|_{L^2(L^2)}^2 \leq$$
$$C \left( \|\mathbf{u}_h\|_{L^2(H^1)}^2 + \|p_h\|_{L^2(J)}^2 + \left\|\mathbf{f}^N\right\|_{0,\Omega}^2 + \|\mathbf{f}_t\|_{L^2(L^2)}^2 + \left\|\mathbf{u}_h^0\right\|_{0,\Omega}^2 + \right.$$
$$\left. |p_h^0|_{J,\Omega}^2 + \left\|\mathbf{f}^1\right\|_{L^2(L^2)}^2 + \|\mathbf{f}\|_{L^2(L^2)}^2 + \|\mathbf{b}\|_{L^2(L^2)}^2 + \|g\|_{L^2(L^2)}^2 \right).$$

Using assumed regularity of the given data to bound the third term and upwards on the righthand side we obtain

$$\left\|\mathbf{u}_h^N\right\|_{1,\Omega}^2 + |p_h^N|_{J,\Omega}^2 + \|\mathbf{z}_h\|_{L^2(L^2)}^2 + \|p_h\|_{L^2(L^2)}^2 \leq C \left( 1 + \|\mathbf{u}_h\|_{L^2(H^1)}^2 + \|p_h\|_{L^2(J)}^2 \right).$$

Upon applying the Gronwall Lemma to the above inequality we obtain the desired result. □

**4.2. Bound on the divergence of the fluid flux.** In order to bound the divergence of the fluid flux we now define the bilinear form $\mathcal{B}_h^n$. We first show how we derive $\mathcal{B}_h^n$ from the fully-discrete weak form (2.7), for which we know that a solution $(\mathbf{u}_h, \mathbf{z}_h, p_h)$ exists for test functions $(\mathbf{v}_h, \mathbf{w}_h, q_h) \in \mathcal{V}_h^X$. Adding (2.7a) and (2.7b), assuming $\mathbf{t}_N = 0$ on $\Gamma_N$, and summing we have

$$\sum_{n=1}^N a(\mathbf{u}_h^n, \mathbf{v}_h) + \sum_{n=1}^N (k^{-1}\mathbf{z}_h^n, \mathbf{w}_h) - \sum_{n=1}^N (p_h^n, \nabla \cdot \mathbf{v}_h) - \sum_{n=1}^N (p_h^n, \nabla \cdot \mathbf{w}_h)$$
$$= \sum_{n=1}^N (\mathbf{f}^n, \mathbf{v}_h) + \sum_{n=1}^N (\mathbf{b}^n, \mathbf{w}_h) \quad \forall (\mathbf{v}_h, \mathbf{w}_h, q_h) \in \mathcal{V}_h^X. \quad (4.6)$$

For the purposes of this proof we now introduce initial conditions for the fluid flux, $\mathbf{z}^0 \in H_{div}(\Omega)$. We also define the projections, into their respective finite element spaces, $\mathbf{z}_h^0 := \pi_h^0 \mathbf{z}^0$ and $p_h^0 := \pi_h^0 p^0$.

Adding (2.7a) and (2.7b), and summing from 0 to $N-1$, we have

$$\sum_{n=1}^N a(\mathbf{u}_h^{n-1}, \mathbf{v}_h) + \sum_{n=1}^N (k^{-1}\mathbf{z}_h^{n-1}, \mathbf{w}_h) - \sum_{n=1}^N (p_h^{n-1}, \nabla \cdot \mathbf{v}_h) - \sum_{n=1}^N (p_h^{n-1}, \nabla \cdot \mathbf{w}_h)$$
$$= \sum_{n=1}^N (\mathbf{f}^{n-1}, \mathbf{v}_h) + \sum_{n=1}^N (\mathbf{b}^{n-1}, \mathbf{w}_h) \quad \forall (\mathbf{v}_h, \mathbf{w}_h, q_h) \in \mathcal{V}_h^X. \quad (4.7)$$

Taking (2.7c), multiplying by $\Delta t$, and summing we have

$$\sum_{n=1}^N \Delta t(\nabla \cdot \mathbf{u}_{\Delta t,h}^n, q_h) + \sum_{n=1}^N \Delta t(\nabla \cdot \mathbf{z}_h^n, q_h) + \sum_{n=1}^N \Delta t J(p_{\Delta t,h}^n, q_h)$$
$$= \sum_{n=1}^N \Delta t(g^n, q_h) \quad \forall (\mathbf{v}_h, \mathbf{w}_h, q_h) \in \mathcal{V}_h^X. \quad (4.8)$$

Now adding (4.6) and (4.8), and subtracting (4.7) we get

$$\sum_{n=1}^{N} \Delta t \mathcal{B}_h^n[(\mathbf{u}_h, \mathbf{z}_h, p_h), (\mathbf{v}_h, \mathbf{w}_h, q_h)]$$
$$= \sum_{n=1}^{N} \Delta t (\mathbf{f}_{\Delta t}^n, \mathbf{v}_h) + \sum_{n=1}^{N} \Delta t (\mathbf{b}_{\Delta t}^n, \mathbf{w}_h) + \sum_{n=1}^{N} \Delta t (g^n, q_h) \ \forall \ (\mathbf{v}_h, \mathbf{w}_h, q_h) \in \mathcal{V}_h^X, \quad (4.9)$$

where

$$\mathcal{B}_h^n[(\mathbf{u}_h, \mathbf{z}_h, p_h), (\mathbf{v}_h, \mathbf{w}_h, q_h)] = a(\mathbf{u}_{\Delta t, h}^n, \mathbf{v}_h) + (k^{-1}\mathbf{z}_{\Delta t, h}^n, \mathbf{w}_h)$$
$$- (p_{\Delta t, h}^n, \nabla \cdot \mathbf{v}_h) - (p_{\Delta t, h}^n, \nabla \cdot \mathbf{w}_h) + (\nabla \cdot \mathbf{u}_{\Delta t, h}^n, q_h) + (\nabla \cdot \mathbf{z}_h^n, q_h) + J(p_{\Delta t, h}^n, q_h). \quad (4.10)$$

With these preliminaries, we may now bound $\mathcal{B}_h^n$ from below.

LEMMA 4.3. *For all* $\beta > \beta^\star > 0$, $(\mathbf{u}_h, \mathbf{z}_h, p_h)$ *satisfies*

$$\sum_{n=1}^{N} \Delta t \ \mathcal{B}_h^n[(\mathbf{u}_h, \mathbf{z}_h, p_h), (\beta \mathbf{u}_{\Delta t, h}^n + \pi_h^1 v_p, \beta \mathbf{z}_h^n, \beta p_{\Delta t, h}^n + \nabla \cdot \mathbf{z}_h^n)] + \left\| \mathbf{z}_h^0 \right\|_{0,\Omega}^2 \geq$$
$$C \left( \left\| \mathbf{u}_{\Delta t, h} \right\|_{L^2(H^1)}^2 + \left\| \mathbf{z}_h^N \right\|_{0,\Omega}^2 + \left\| p_{\Delta t, h} \right\|_{L^2(L^2)}^2 + \left\| p_{\Delta t, h} \right\|_{L^2(J)}^2 + \left\| \nabla \cdot \mathbf{z}_h \right\|_{L^2(L^2)}^2 \right).$$

*Proof.* For $n = 1, 2, \ldots, N$ we choose $(\mathbf{v}_h, \mathbf{w}_h, q_h) = (\beta \mathbf{u}_{\Delta t,h}^n + \pi_h^1 \mathbf{v}_{p_h^n}, \beta \mathbf{z}_h^n, \beta p_{\Delta t,h}^n + \nabla \cdot \mathbf{z}_h^n)$ in (4.10)

$$\sum_{n=1}^{N} \Delta t \mathcal{B}_h^n[(\mathbf{u}_h, \mathbf{z}_h, p_h), (\beta \mathbf{u}_{\Delta t, h}^n + \pi_h^1 \mathbf{v}_p, \beta \mathbf{z}_h^n, \beta p_{\Delta t, h}^n + \nabla \cdot \mathbf{z}_h^n)]$$
$$= \sum_{n=1}^{N} \Delta t a(\mathbf{u}_{\Delta t, h}^n, \beta \mathbf{u}_{\Delta t, h}^n) + \sum_{n=1}^{N} \Delta t (k^{-1}\mathbf{z}_{\Delta t, h}^n, \beta \mathbf{z}_h^n) + \sum_{n=1}^{N} \Delta t (\nabla \cdot \mathbf{z}_h^n, \nabla \cdot \mathbf{z}_h^n)$$
$$+ \sum_{n=1}^{N} \Delta t (\mathbf{u}_{\Delta t, h}^n, \nabla \cdot \mathbf{z}_h^n) + \sum_{n=1}^{N} \Delta t J(p_{\Delta t, h}^n, \nabla \cdot \mathbf{z}_h^n) + \sum_{n=1}^{N} \Delta t J(p_{\Delta t, h}^n, \beta p_{\Delta t, h}^n)$$
$$+ \sum_{n=1}^{N} \Delta t a(\mathbf{u}_{\Delta t, h}^n, \pi_h^1 \mathbf{v}_p) - \sum_{n=1}^{N} \Delta t (p_{\Delta t, h}^n, \nabla \cdot \pi_h^1 \mathbf{v}_p). \quad (4.11)$$

For all $\epsilon > 0$ using (2.3), (2.4), the Cauchy-Schwarz, Young's and Poincaré inequalities, (3.1) on $\nabla \cdot \mathbf{z}_h^n$, and an approach similar to step 2 in the proof of Theorem 3.2 for the final two terms on the righthand side, we obtain

$$\sum_{n=1}^{N} \Delta t \mathcal{B}_h^n[(\mathbf{u}_h, \mathbf{z}_h, p_h), (\beta \mathbf{u}_{\Delta t, h}^n + \pi_h^1 \mathbf{v}_p, \beta \mathbf{z}_h^n, \beta p_{\Delta t, h}^n + \nabla \cdot \mathbf{z}_h^n)]$$
$$\geq \left( \beta C_k - \frac{C_p + C_c}{2\epsilon} \right) \left\| \mathbf{u}_{\Delta t, h} \right\|_{L^2(H^1)}^2 + \frac{\beta \lambda_{max}^{-1}}{2} \left\| \mathbf{z}_h^N \right\|_{0,\Omega}^2 + \left( \beta - \frac{3}{4\epsilon} \right) \left\| p_{\Delta t, h} \right\|_{L^2(J)}^2$$
$$+ (1 - \epsilon(1 + c_z)) \left\| \nabla \cdot \mathbf{z}_h \right\|_{L^2(L^2)}^2 - \frac{\beta \lambda_{min}^{-1}}{2} \left\| \mathbf{z}_h^0 \right\|_{0,\Omega}^2 + (1 - C\epsilon) \left\| p_{\Delta t, h} \right\|_{L^2(L^2)}^2. \quad (4.12)$$

Finally choosing $\epsilon$ sufficiently small and $\beta \geq \max\left[\frac{C_p}{2C_k\epsilon}, \frac{3}{4\epsilon}\right]$ completes the proof. □

The following Lemma shows the divergence control of the fluid flux.

LEMMA 4.4. $\mathbf{z}_h$ obtained from (4.9) satisfies

$$\|\nabla \cdot \mathbf{z}_h\|^2_{L^2(L^2)} \leq C.$$

*Proof.* For $n = 1, 2, \ldots, N$ we choose $(\mathbf{v}_h, \mathbf{w}_h, q_h) = (\beta \mathbf{u}^n_{\Delta t,h} + \pi^1_h \mathbf{v}_{p^n_h}, \beta \mathbf{z}^n_h, \beta p^n_{\Delta t,h} + \nabla \cdot \mathbf{z}^n_h)$ in (4.9) yielding

$$\sum_{n=1}^N \Delta t \mathcal{B}^n_h[(\mathbf{u}^n_h, \mathbf{z}^n_h, p^n_h), (\beta \mathbf{u}^n_{\Delta t,h} + \pi^1_h \mathbf{v}_{p^n_h}, \mathbf{z}^n_h, \beta p^n_{\Delta t,h} + \nabla \cdot \mathbf{z}^n_h)]$$
$$= \sum_{n=1}^N \Delta t (\mathbf{f}^n_{\Delta t}, \beta \mathbf{u}^n_{\Delta t,h} + \pi^1_h \mathbf{v}_{p^n_h}) + \sum_{n=1}^N \Delta t (\mathbf{b}^n_{\Delta t}, \beta \mathbf{z}^n_h) + \sum_{n=1}^N \Delta t (g^n, \beta p^n_{\Delta t,h} + \nabla \cdot \mathbf{z}^n_h).$$

Using Lemma 4.3, the Cauchy-Schwarz and Young's inequalities, and (3.6), along with ideas already presented in the proof of Lemma 4.2

$$\|\mathbf{u}_{\Delta t,h}\|^2_{L^2(H^1)} + \|p_{\Delta t,h}\|^2_{L^2(L^2)} + \|p_{\Delta t,h}\|^2_{L^2(J)} + \left\|\mathbf{z}^N_h\right\|^2_{0,\Omega} + \|\nabla \cdot \mathbf{z}_h\|^2_{L^2(L^2)}$$
$$\leq C \left(\|\mathbf{f}_t\|^2_{L^2(L^2)} + \|\mathbf{b}_t\|^2_{L^2(L^2)} + \|p_h\|^2_{L^2(L^2)} + \|\mathbf{z}_h\|^2_{L^2(L^2)} + \|g\|^2_{L^2(L^2)}\right).$$

Finally, by applying a Gronwall Lemma, using Lemma 4.2 and regularity, we obtain the desired result. □

**4.3. The energy estimate.** THEOREM 4.5. *The solution to the fully-discrete problem (2.7) satisfies the energy estimate*

$$\|\mathbf{u}_h\|^2_{L^\infty(H^1)} + \|p_h\|^2_{L^\infty(J)} + \|\mathbf{z}_h\|^2_{L^2(L^2)} + \|p_h\|^2_{L^2(L^2)} + \|\nabla \cdot \mathbf{z}_h\|^2_{L^2(L^2)} \leq C.$$

*Proof.* The proof follows from combining Lemma 4.2 and Lemma 4.4, and noting that these Lemmas hold for all time steps $n = 0, 1, ..., N$. This then gives the desired discrete in time $L^\infty$ bounds. □

REMARK 4.1. *Having proven Theorem 4.5, it is now a standard calculation to show that the discrete Galerkin approximation converges weakly, as $\Delta t, h \to 0$, to the continuous problem with respect to continuous versions of the norms of the energy estimate in Theorem 4.5. This in turn shows that the continuous variational problem is well-posed. Due to the linearity of the variational form and noting that $|\mathbf{v}|_{J,\Omega} \to 0$ as $h \to 0$, these calculations are straight forward and closely follow the existence and uniqueness proofs presented in [39] and [4] for the linear two-field Biot problem and a nonlinear Biot problem, respectively.*

**5. A-priori error analysis.** Lemma 5.1 provides a Galerkin orthogonality result obtained by comparing continuous and discrete weak forms, which is the corner stone of the error analysis. Lemma 5.2 bounds the auxiliary errors for displacement, flux and pressure in the appropriate norms and Lemma 5.3 bounds the auxiliary error for the divergence of the flux. Since Lemmas 5.2 and 5.3 bound the auxiliary errors at the same order as the projection errors, combining projection and auxiliary errors in Theorem 5.4 provides an optimal error estimate.

We define the finite element error functions

$$\mathbf{e_u} := \mathbf{u} - \mathbf{u}_h, \quad \mathbf{e_z} := \mathbf{z} - \mathbf{z}_h, \quad e_p := p - p_h.$$

We introduce the following projection errors:

$$\eta_{\mathbf{u}} := \mathbf{u} - \pi_h^1 \mathbf{u}, \ \eta_{\mathbf{z}} := \mathbf{z} - \pi_h^1 \mathbf{z}, \ \eta_p := p - \pi_h^0 p,$$

where we have assumed $\mathbf{z}(t_n, \cdot) \in (H^1(\Omega))^d$.

Auxiliary errors:

$$\theta_{\mathbf{u}}^n(\cdot) := \pi_h^1 \mathbf{u}(t_n, \cdot) - \mathbf{u}_h^n(\cdot), \ \theta_{\mathbf{z}}^n(\cdot) := \pi_h^1 \mathbf{z}(t_n, \cdot) - \mathbf{z}_h^n(\cdot), \ \theta_p^n(\cdot) := \pi_h^0 p(t_n, \cdot) - p_h^n(\cdot), \tag{5.1}$$

and time-discretization errors:

$$\rho_{\mathbf{u}}^n(\cdot) := \frac{\mathbf{u}(t_n, \cdot) - \mathbf{u}(t_{n-1}, \cdot)}{\Delta t} - \frac{\partial \mathbf{u}(t_n, \cdot)}{\partial t}, \ \rho_p^n := \frac{p(t_n, \cdot) - p(t_{n-1}, \cdot)}{\Delta t} - \frac{\partial p(t_n, \cdot)}{\partial t}. \tag{5.2}$$

**5.1. Galerkin orthogonality.** We now give a Galerkin orthogonality type argument for analysing the difference between the fully-discrete approximation and the true solution. For this we introduce the continuous counterpart of the fully-discrete combined weak form (4.1) given by

$$B^n[(\mathbf{u}, \mathbf{z}, p), (\mathbf{v}, \mathbf{w}, q)] = (\mathbf{f}(t_n, \cdot), \mathbf{v}) + (\mathbf{b}(t_n, \cdot), \mathbf{w}) + (g(t_n, \cdot), q) \ \forall \ (\mathbf{v}, \mathbf{w}, q) \in \mathcal{V}^X, \tag{5.3}$$

where

$$\begin{aligned} B^n[(\mathbf{u}, \mathbf{z}, p), (\mathbf{v}, \mathbf{w}, q)] = &\ a(\mathbf{u}(t_n, \cdot), \mathbf{v}) + (k^{-1} \mathbf{z}(t_n, \cdot), \mathbf{w}) - (p(t_n, \cdot), \nabla \cdot \mathbf{v}) \\ &- (p(t_n, \cdot), \nabla \cdot \mathbf{w}) + (\nabla \cdot \mathbf{u}_t(t_n, \cdot), q) + (\nabla \cdot \mathbf{z}(t_n, \cdot), q). \end{aligned}$$

LEMMA 5.1. *Assuming*
$(\mathbf{u}(t_n, \cdot), \mathbf{z}(t_n, \cdot), p(t_n, \cdot)) \in \left(H^1(\Omega)\right)^d \times H_{div}(\Omega) \times \left(H^1(\Omega) \cap \mathcal{L}(\Omega)\right),$

$$B_{\Delta t, h}^n[(\mathbf{e_u}, \mathbf{e_z}, e_p), (\mathbf{v}_h, \mathbf{w}_h, q_h)] = (\nabla \cdot \rho_{\mathbf{u}}^n, q_h) + J(\rho_p^n, q_h) \ \ \forall (\mathbf{v}_h, \mathbf{w}_h, q_h) \in \mathcal{V}_h^X.$$

*Proof.* Subtracting the discrete weak form (4.1) from the continuous weak form (5.3), we obtain

$$B^n[(\mathbf{u}, \mathbf{z}, p), (\mathbf{v}_h, \mathbf{w}_h, q_h)] - B_{\Delta t, h}^n[(\mathbf{u}_h, \mathbf{z}_h, p_h), (\mathbf{v}_h, \mathbf{w}_h, q_h)] = 0, \ \ \forall (\mathbf{v}_h, \mathbf{w}_h, q_h) \in \mathcal{V}_h^X.$$

Adding $J(p_t(t_n, \cdot), q) = 0$ to the left, see (3.2) and $(\nabla \cdot (\mathbf{u}_{\Delta t}(t_n, \cdot) - \mathbf{u}_t(t_n, \cdot)), q) + J(p_{\Delta t}(t_n, \cdot) - p_t(t_n, \cdot), q)$ to both the left and righthand sides we obtain the desired result. □

**5.2. Auxiliary error estimates.** LEMMA 5.2. *Assuming*
$\mathbf{u} \in H^2\left(0, T; \left(L^2(\Omega)\right)^d\right) \cap H^1\left(0, T; \left(H^2(\Omega)\right)^d\right), \ \mathbf{z} \in L^2\left(0, T; \left(H^1(\Omega)\right)^d\right)$ and $p \in H^2\left(0, T; H^1(\Omega) \cap \mathcal{L}(\Omega)\right),$ *then the finite element solution (2.7) satisfies the error estimate*

$$\interleave[\theta_{\mathbf{u}}, \theta_{\mathbf{z}}, \theta_p]\interleave_B^2 + \|\theta_p\|_{L^\infty(J)}^2 \leq C(T)(h^2 + \Delta t^2). \tag{5.4}$$

*Proof.* Using Lemma 5.1 and choosing $\mathbf{v}_h^n = \theta_{\Delta t, \mathbf{u}}^n + \pi_h^1 \mathbf{v}_{p_h^n}$, $\mathbf{w}_h^n = \theta_\mathbf{z}^n$, $q_h^n = \theta_p^n$, we get

$$B_{\Delta t, h}^n [(\theta_\mathbf{u}^n + \eta_\mathbf{u}^n, \theta_\mathbf{z}^n + \eta_\mathbf{z}^n, \theta_p^n + \eta_p^n), (\theta_{\Delta t, \mathbf{u}}^n + \pi_h^1 \mathbf{v}_{p_h^n}, \theta_\mathbf{z}^n, \theta_p^n)]$$
$$= (\nabla \cdot \rho_\mathbf{u}^n, \theta_p^n) + J(\rho_p^n, \theta_p^n).$$

Rearranging gives

$$B_{\Delta t, h}^n [(\theta_\mathbf{u}^n, \theta_\mathbf{z}^n, \theta_p^n), (\theta_{\Delta t, \mathbf{u}}^n + \pi_h^1 \mathbf{v}_{p_h^n}, \theta_\mathbf{z}^n, \theta_p^n)]$$
$$= (\nabla \cdot \rho_\mathbf{u}^n, \theta_p^n) + J(\rho_p^n, \theta_p^n) - B_{\Delta t, h}^n [(\eta_\mathbf{u}^n, \eta_\mathbf{z}^n, \eta_p^n), (\theta_{\Delta t, \mathbf{u}}^n + \pi_h^1 \mathbf{v}_{p_h^n}, \theta_\mathbf{z}^n, \theta_p^n)].$$

Expanding the righthand side, noting that $(\eta_p^n, \nabla \cdot (\theta_{\Delta t, \mathbf{u}}^n + \pi_h^1 \mathbf{v}_p)) = 0$, $(\eta_p^n, \nabla \cdot \theta_\mathbf{z}^n) = 0$, multiplying both sides by $\Delta t$ and summing gives

$$\sum_{n=1}^N \Delta t B_{\Delta t, h}^n [(\theta_\mathbf{u}^n, \theta_\mathbf{z}^n, \theta_p^n), (\theta_{\Delta t, \mathbf{u}}^n + \pi_h^1 \mathbf{v}_{p_h^n}, \theta_\mathbf{z}^n, \theta_p^n)] = \sum_{i=1}^7 \Phi_i,$$

where

$$\Phi_1 := -\sum_{n=1}^N \Delta t a(\eta_\mathbf{u}^n, \theta_{\Delta t, \mathbf{u}}^n), \quad \Phi_2 := -\sum_{n=1}^N \Delta t ((k^{-1} \eta_\mathbf{z}^n, \theta_\mathbf{z}^n)), \quad \Phi_3 := -\sum_{n=1}^N \Delta t a(\eta_\mathbf{u}^n, \pi_h^1 \mathbf{v}_p),$$

$$\Phi_4 := -\sum_{n=1}^N \Delta t J(\eta_{\Delta t, p}^n, \theta_p^n), \quad \Phi_5 := \sum_{n=1}^N \Delta t (\nabla \cdot \rho_\mathbf{u}^n, \theta_p^n), \quad \Phi_6 := \sum_{n=1}^N \Delta t J(\rho_p^n, \theta_p^n),$$

$$\Phi_7 := -\sum_{n=1}^N \Delta t (\theta_p^n, \nabla \cdot (\eta_{\Delta t, \mathbf{u}}^n + \eta_\mathbf{z}^n)).$$

We now individually consider the terms on the right hand side of (5.5):

To bound the first quantity, we use (3.11), Lemma 3.1, the triangle, Cauchy-Schwarz and Young's inequalities, $\theta_\mathbf{u}^0 = 0$, and (2.2),

$$\Phi_1 = -\sum_{n=1}^N a(\eta_\mathbf{u}^n, \theta_\mathbf{u}^n - \theta_\mathbf{u}^{n-1})$$
$$= -a(\eta_\mathbf{u}^N, \theta_\mathbf{u}^N) + \sum_{n=1}^N a(\eta_\mathbf{u}^n - \eta_\mathbf{u}^{n-1}, \theta_\mathbf{u}^{n-1})$$
$$= -a(\eta_\mathbf{u}^N, \theta_\mathbf{u}^N) + \Delta t \sum_{n=1}^N a\left((I - \pi_h^1)\left(\rho_\mathbf{u}^n + \frac{\partial \mathbf{u}(t_n, \cdot)}{\partial t}\right), \theta_\mathbf{u}^{n-1}\right)$$
$$\leq \epsilon C \|\theta_\mathbf{u}^N\|_{1, \Omega}^2 + \frac{Ch^2}{\epsilon} \|\mathbf{u}^N\|_{2, \Omega}^2 + \epsilon C \|\theta_\mathbf{u}\|_{L^2(H^1)}^2$$
$$+ \frac{Ch^2}{2\epsilon} \|\mathbf{u}_t\|_{L^2(H^2)}^2 + \frac{C \Delta t^2}{2\epsilon} \|\mathbf{u}_{tt}\|_{L^2(H^1)}^2.$$
(5.5)

Next, using (2.4), Young's inequality, (3.6) and Lemma 3.1,

$$\Phi_2 \leq \frac{\epsilon}{2} \|\theta_\mathbf{z}\|_{L^2(L^2)}^2 + \frac{\lambda_{min}^{-2} h^2}{2\epsilon} \|\mathbf{z}\|_{L^2(H^1)}^2.$$

Using (2.2), Young's inequality and Lemma 3.1,

$$\Phi_3 \leq \frac{\epsilon}{2}\|\pi_h^1 \mathbf{v}_{p_h^n}\|^2_{L^2(H^1)} + \frac{C}{2\epsilon}\|\eta_\mathbf{u}\|^2_{L^2(H^1)} \leq \frac{\epsilon \hat{c}^2}{2}\|\theta_p\|^2_{L^2(L^2)} + \frac{Ch^2}{2\epsilon}\|\mathbf{u}\|^2_{L^2(H^2)}.$$

The bound on $\Phi_4$ is obtained using a similar argument to the bound on $\Phi_1$,

$$\Phi_4 \leq \epsilon\|\theta_p\|^2_{L^2(J)} + \frac{h^2}{2\epsilon}\|p_t\|^2_{L^2(H^1)} + \frac{\Delta t^2}{2\epsilon}\|p_{tt}\|^2_{L^2(H^1)}.$$

Using the Cauchy-Schwarz and Young's inequalities and Lemma 3.1,

$$\Phi_5 \leq \frac{\epsilon}{2}\|\theta_p\|^2_{L^2(L^2)} + \frac{\Delta t^2}{2\epsilon}\|\mathbf{u}_{tt}\|^2_{L^2(L^2)} \text{ and } \Phi_6 \leq \frac{\epsilon}{2}\|\theta_p\|^2_{L^2(J)} + \frac{\Delta t^2}{2\epsilon}\|p_{tt}\|^2_{L^2(L^2)}.$$

Finally, using the Cauchy-Schwarz and Young's inequalities, and a similar argument to the bound on $\Phi_1$,

$$\Phi_7 \leq \frac{3\epsilon}{2}\|\theta_p\|^2_{L^2(L^2)} + \frac{h^2}{2\epsilon}\|\mathbf{u}_t\|^2_{L^2(H^2)} + \frac{\Delta t^2}{2\epsilon}\|\mathbf{u}_{tt}\|^2_{L^2(H^1)} + \frac{h^2}{2\epsilon}\|\mathbf{z}\|^2_{L^2(H^2)}.$$

Combining these bounds with an application of the coercivity Lemma 4.1 to (5.5), noting the assumed regularity of the continuous solution and choosing $\epsilon$ sufficiently small, gives

$$\|\theta_\mathbf{u}^N\|^2_{1,\Omega} + |\theta_p^N|^2_{J,\Omega} + \|\theta_\mathbf{z}\|^2_{L^2(L^2)} + \|\theta_p\|^2_{L^2(L^2)} \leq C\left(\|\theta_\mathbf{u}\|^2_{L^2(H^1)} + \|\theta_p\|^2_{L^2(J)} + h^2 + \Delta t^2\right). \tag{5.6}$$

An application of Gronwall's Lemma gives

$$\|\theta_\mathbf{u}^N\|^2_{1,\Omega} + |\theta_p^N|^2_{J,\Omega} + \|\theta_\mathbf{z}\|^2_{L^2(L^2)} + \|\theta_p\|^2_{L^2(L^2)} \leq C(T)\left(h^2 + \Delta t^2\right).$$

Because the above holds for all time steps $n = 0, 1, ..., N$, we can get the desired $L^\infty$ bounds to complete the proof of the theorem. □

We now present an a-priori auxiliary error estimate of the fluid flux, in its natural $Hdiv$ norm.

LEMMA 5.3. *Assuming* $\mathbf{u} \in H^2\left(0, T; \left(H^1(\Omega)\right)^d\right) \cap H^1\left(0, T; \left(H^2(\Omega)\right)^d\right)$, $\mathbf{z} \in L^2\left(0, T; \left(H^2(\Omega)\right)^d\right)$ *and* $p \in H^2\left(0, T; J \cap \mathcal{L}(\Omega)\right) \cap H^1(0, T; H^1(\Omega))$, *then the finite element solution (2.7) satisfies the auxillary error estimate*

$$\|\nabla \cdot \theta_\mathbf{z}\|^2_{L^2(L^2)} \leq C(T)(h^2 + \Delta t^2). \tag{5.7}$$

*Proof.*
Similarly to the approach taken in obtaining (4.9) we may easily obtain the following identity

$$\sum_{n=1}^{N} \Delta t \mathcal{B}_h^n[(\theta_\mathbf{u}^n, \theta_\mathbf{z}^n, \theta_p^n), (\beta\theta_{\Delta t,\mathbf{u}}^n + \pi_h^1 \mathbf{v}_{\theta_{\Delta t,p}^n}, \beta\theta_\mathbf{z}^n, \beta\theta_{\Delta t,p}^n + \nabla \cdot \theta_\mathbf{z}^n)] = \sum_{i=1}^{6} \Psi_i,$$

where

$$\Psi_1 := -\sum_{n=1}^{N} \Delta t a(\eta^n_{\Delta t,\mathbf{u}}, \beta\theta^n_{\Delta t,\mathbf{u}} + \pi^1_h \mathbf{v}_{\theta^n_{\Delta t,p}}),$$

$$\Psi_2 := -\sum_{n=1}^{N} \Delta t(\nabla \cdot (\eta^n_{\Delta t,\mathbf{u}} + \eta^n_{\mathbf{z}}), \nabla \cdot \theta^n_{\mathbf{z}} + \beta\theta^n_{\Delta t,p}),$$

$$\Psi_3 := \sum_{n=1}^{N} \Delta t J(\eta^n_{\Delta t,p}, \beta\theta^n_{\Delta t,p} + \nabla \cdot \theta^n_{\mathbf{z}}), \quad \Psi_4 := -\sum_{n=1}^{N} \Delta t((k^{-1}\eta^n_{\Delta t,\mathbf{z}}, \beta\theta^n_{\mathbf{z}})),$$

$$\Psi_5 := \sum_{n=1}^{N} \Delta t J(\rho^n_p, \beta\theta^n_{\Delta t,p} + \nabla \cdot \theta^n_{\mathbf{z}}), \quad \Psi_6 := \sum_{n=1}^{N} \Delta t(\nabla \cdot \rho^n_{\mathbf{u}}, \beta\theta^n_{\Delta t,p} + \nabla \cdot \theta^n_{\mathbf{z}}).$$

We now bound the terms on the right hand side of (5.8) using machinery developed during the previous proof:

$$\Psi_1 \leq \frac{C\epsilon}{2}\|\theta_{\Delta t,\mathbf{u}}\|^2_{L^2(H^1)} + \frac{\hat{c}^2\epsilon}{2}\|\theta_{\Delta t,p}\|^2_{L^2(L^2)} + \frac{Ch^2}{2\epsilon}\|\mathbf{u}_t\|^2_{L^2(H^2)}$$
$$+ \frac{C}{2\epsilon}\Delta t^2 \|\mathbf{u}_{tt}\|^2_{L^2(H^1)}, \tag{5.8}$$

$$\Psi_2 \leq \epsilon\|\nabla \cdot \theta_{\mathbf{z}}\|^2_{L^2(L^2)} + \epsilon\|\theta_{\Delta t,p}\|^2_{L^2(L^2)} + \frac{Ch^2}{2\epsilon}\left(\|\mathbf{u}_t\|^2_{L^2(H^2)} + \|\mathbf{z}\|^2_{L^2(H^2)}\right)$$
$$+ \frac{C}{2\epsilon}\Delta t^2 \|\mathbf{u}_{tt}\|^2_{L^2(H^1)}, \tag{5.9}$$

$$\Psi_3 \leq \epsilon C\|\nabla \cdot \theta_{\mathbf{z}}\|^2_{L^2(L^2)} + \epsilon\|\theta^n_{\Delta t,p}\|^2_{L^2(J)} + \frac{Ch^2}{2\epsilon}\|p_t\|^2_{L^2(H^1)} + \frac{C}{2\epsilon}\Delta t^2\|p_{tt}\|^2_{L^2(J)} \tag{5.10}$$

$$\Psi_4 \leq \epsilon\|\theta_{\mathbf{z}}\|^2_{L^2(L^2)} + \frac{Ch^2}{2\epsilon}\|\mathbf{z}_t\|^2_{L^2(H^1)} + \frac{C}{2\epsilon}\Delta t^2\|\mathbf{z}_{tt}\|^2_{L^2(L^2)}, \tag{5.11}$$

$$\Psi_5 \leq \epsilon\|\theta_{\Delta t,p}\|^2_{L^2(J)} + \epsilon C\|\nabla \cdot \theta_{\mathbf{z}}\|^2_{L^2(L^2)} + \frac{C\Delta t^2}{2\epsilon}\|p_{tt}\|^2_{L^2(J)}, \tag{5.12}$$

$$\Psi_6 \leq \epsilon\|\theta_{\Delta t,p}\|^2_{L^2(L^2)} + \epsilon\|\nabla \cdot \theta_{\mathbf{z}}\|^2_{L^2(L^2)} + \frac{C}{2\epsilon}\Delta t^2\|\mathbf{u}_{tt}\|^2_{L^2(H^1)}. \tag{5.13}$$

We can now combine the individual bounds (5.8), (5.9), (5.10), (5.11), (5.12), and (5.13), with the coercivity result Lemma 4.3, choose $\beta$ sufficiently large, use the assumption $\theta^0_{\mathbf{z}} = 0$, the assumed regularity of $\mathbf{u}, \mathbf{z}$ and $p$, and choose $\epsilon$ sufficiently small to obtain

$$\left\|\theta^N_{\mathbf{z}}\right\|^2_{0,\Omega} + \|\nabla \cdot \theta_{\mathbf{z}}\|^2_{L^2(L^2)} \leq C\|\theta_{\mathbf{z}}\|^2_{L^2(L^2)} + C(h^2 + \Delta t^2).$$

Applying Gronwall's Lemma, we get the desired result. □

**5.3. The a priori error estimate.** By combining the previous Lemmas we obtain the main Theorem regarding the finite element error estimate.

THEOREM 5.4. *If* $\mathbf{u} \in H^2\left(0,T;\left(L^2(\Omega)\right)^d\right) \cap H^1\left(0,T;\left(H^2(\Omega)\right)^d\right)$, $\mathbf{z} \in L^2\left(0,T;\left(H^1(\Omega)\right)^d\right)$ *and* $p \in H^2\left(0,T; H^1(\Omega) \cap \mathcal{L}(\Omega)\right)$, *then finite element solution (2.7) satisfies the error estimate*

$$\|\!|\mathbf{e}_\mathbf{u}, \mathbf{e}_\mathbf{z}, e_p|\!\|^2_B \leq C(h^2 + \Delta t^2).$$

*Moreover, if* $\mathbf{u} \in H^2\left(0, T; \left(H^1(\Omega)\right)^d\right) \cap H^1\left(0, T; \left(H^2(\Omega)\right)^d\right)$, $\mathbf{z} \in L^2\left(0, T; \left(H^2(\Omega)\right)^d\right)$ *and* $p \in H^2\left(0, T; J \cap \mathcal{L}(\Omega)\right) \cap H^1(0, T; H^1(\Omega))$, *then the finite element solution (2.7) satisfies the error estimate*

$$\||\mathbf{e_u}, \mathbf{e_z}, e_p\||_B^2 + \|\nabla \cdot \mathbf{e_z}\|_{L^2(L^2)}^2 \leq C(h^2 + \Delta t^2).$$

*Proof.* We first write the errors as $\mathbf{e_u}^n = \eta_\mathbf{u}^n + \theta_\mathbf{u}^n$, and similarly for the other variables. Using Lemma 3.1 we can bound the projection errors, and using Lemma 5.2 and Lemma 5.3 we can bound the auxillary errors to give the desired result. □

**6. Numerical Results.** We first present convergence studies for both two- and three-dimensional test problems which illustrate the predicted convergence rates for the fully-discrete finite element method. We then apply our method to the popular 2D cantilever bracket problem and demonstrate that our stabilization technique overcomes the spurious pressure oscillations that have been experienced by other methods. Finally, a 3D unconfined compression problem is presented that highlights the added mass effect of the method for different choices of the stabilization parameter $\delta$.

**6.1. 2D test problem.** Choosing $\lambda = \mu = \alpha = 1$ and $c_0 = 0$ in (2.1) we solve the problem

$$-2\nabla(\nabla \cdot \mathbf{u}) - \nabla^2 \mathbf{u} + \nabla p = \mathbf{f} \quad \text{in } \Omega, \tag{6.1a}$$

$$\mathbf{z} + \nabla p = 0 \quad \text{in } \Omega, \tag{6.1b}$$

$$\nabla \cdot (\mathbf{u}_t + \mathbf{z}) = g \quad \text{in } \Omega, \tag{6.1c}$$

$$\mathbf{u}(t) = \mathbf{u}_D \quad \text{on } \Gamma_D = \Gamma, \tag{6.1d}$$

$$\mathbf{z}(t) \cdot \mathbf{n} = q_D \quad \text{on } \Gamma_F = \Gamma, \tag{6.1e}$$

$$\mathbf{u}(0, \mathbf{x}) = 0, \quad p(0, \mathbf{x}) = 0 \quad \mathbf{x} \in \Omega. \tag{6.1f}$$

The domain, $\Omega$, is the unit square and the source terms and boundary conditions are chosen so that the true solution is

$$\mathbf{u} = \begin{pmatrix} -\frac{1}{4\pi} \cos(2\pi x) \sin(2\pi y) \sin(2\pi t) \\ -\frac{1}{4\pi} \sin(2\pi x) \cos(2\pi y) \sin(2\pi t) \end{pmatrix}, \quad \mathbf{z} = \begin{pmatrix} -2\pi \cos(2\pi x) \sin(2\pi y) \sin(2\pi t) \\ -2\pi \sin(2\pi x) \cos(2\pi y) \sin(2\pi t) \end{pmatrix},$$

and $p = \sin(2\pi x) \sin(2\pi y) \sin(2\pi t)$, with $t \in [0, 0.25]$.

**6.1.1. Choice of $\delta$.** The most appropriate choice of stabilization parameter $\delta$ is not known a priori. Small values of $\delta$ can result in spurious pressure solutions, as shown in Figure 6.1a for $\delta = 0.1$. Larger values of the stabilization parameter produce smooth pressure solutions, as shown in Figure 6.1b for a value of $\delta = 1$. The value of $\delta$ required to produce a stable solution depends on the geometry and material parameters of the particular problem under investigation, but is independent of any mesh parameters.

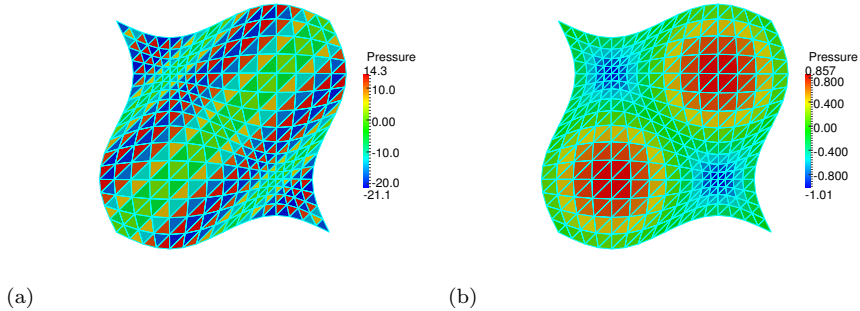

(a)        (b)

Fig. 6.1: (a) Unstable pressure field with $\delta = 0.1$ at $t = 0.25$, stabilization parameter too small. (b) Stable pressure field, with $\delta = 1$ at $t = 0.25$.

**6.1.2. 2D convergence study.** The convergence of the method with discretization parameters is illustrated in Figure 6.2a – 6.2e for $\delta = 1, 10, 100$. The convergence rates observed in the appropriate norms agree with the theoretically derived error estimates.

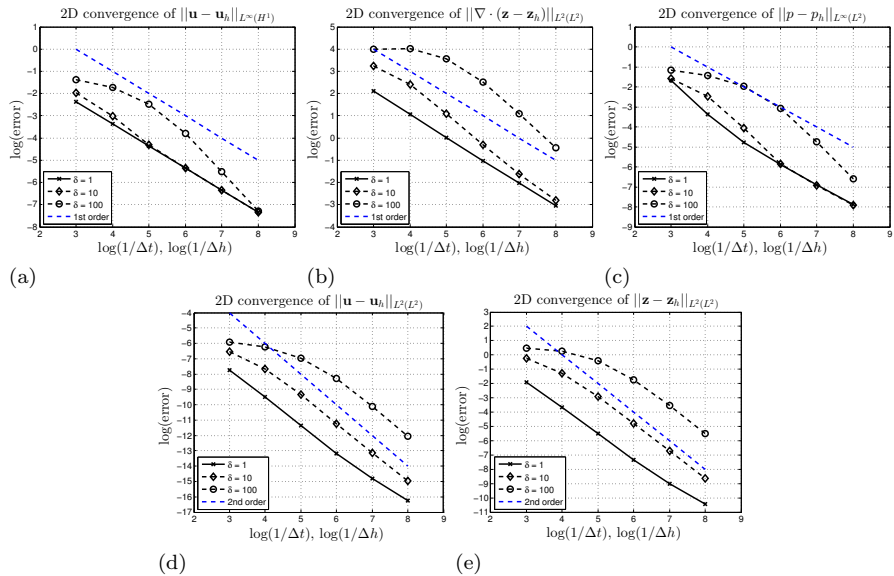

Fig. 6.2: Convergence of the displacement, fluid flux, and pressure errors in their respective norms for the simplified poroelastic 2D test problem with different (stable) values for the stabilization parameter $\delta$.

**6.1.3. Alternative stabilization techniques.** In Figure 6.3 we illustrate the convergence of the pressure error for three possible stabilization forms. As demonstrated in Section 6.1, the stabilization $J(p_{\Delta t,h}, q_h)$ yields a stable solution and optimal convergence rate. A more naive approach, inserting the stabilization $J(p_h, q_h)$ results in the solution becoming unstable after the first refinement step. This is because the stabilization becomes relatively small as $\Delta t$ decreases. To overcome this issue one could chose to scale the stabilization, and try $\frac{1}{\Delta t} J(p_h, q_h)$. Although this stabilization now stays stable during refinement, it does not converge at an optimal rate.

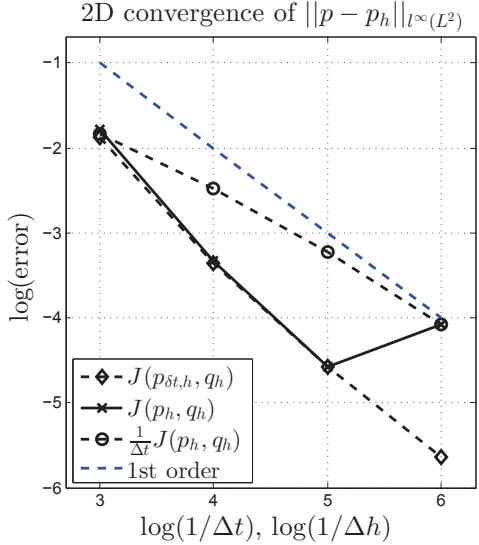

Fig. 6.3: Convergence of the pressure error for three different stabilization forms, with $\delta = 1$.

**6.2. 3D test problem.** Extending the test problem in Section 6.1 to the unit cube, we set

$$\mathbf{u} = \begin{pmatrix} -\frac{1}{6\pi} \cos(2\pi x) \sin(2\pi y) \sin(2\pi z) \sin(2\pi t) \\ -\frac{1}{6\pi} \sin(2\pi x) \cos(2\pi y) \sin(2\pi z) \sin(2\pi t) \\ -\frac{1}{6\pi} \sin(2\pi x) \sin(2\pi y) \cos(2\pi z) \sin(2\pi t) \end{pmatrix},$$

$$\mathbf{z} = \begin{pmatrix} -2\pi \cos(2\pi x) \sin(2\pi y) \sin(2\pi z) \sin(2\pi t) \\ -2\pi \sin(2\pi x) \cos(2\pi y) \sin(2\pi z) \sin(2\pi t) \\ -2\pi \sin(2\pi x) \sin(2\pi y) \cos(2\pi z) \sin(2\pi t) \end{pmatrix},$$

and

$$p = \sin(2\pi x) \sin(2\pi y) \sin(2\pi z) \sin(2\pi t).$$

The expected rates of convergence for each variable in the appropriate norm are illustrated in the numerical results presented in Figure 6.4a – 6.4e for $\delta = 0.001, 0.01, 0.1$. The stabilization factor $\delta$ may be chosen to be very much smaller for 3D problems as compared to 2D problems and the effect of the stabilization term on the solution is negligible. This can be explained by the improved ratio of solid displacement and fluid flux nodes to pressure nodes in three dimensions, making the LBB condition easier to satisfy.

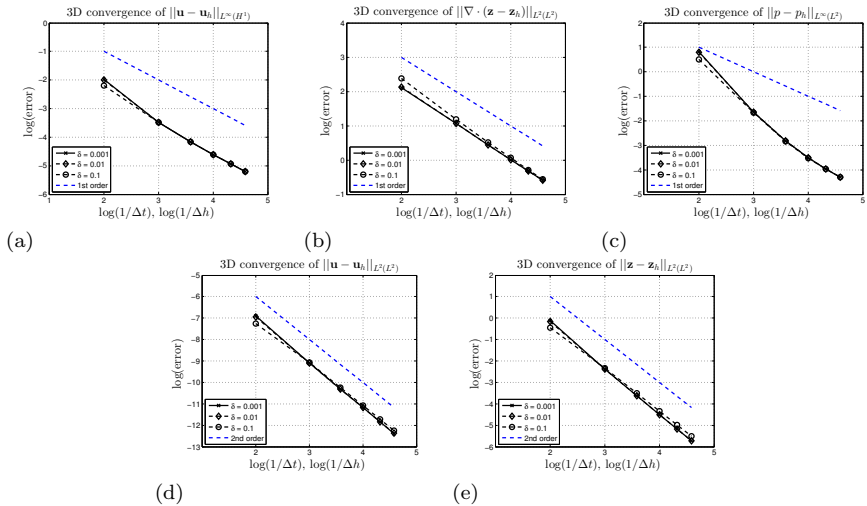

Fig. 6.4: Convergence of the displacement, fluid flux, and pressure errors in their respective norms for the simplified poroelastic 3D test problem with different (stable) values for the stabilization parameter $\delta$.

**6.3. 2D cantilever bracket problem.** We consider the 2D cantilever bracket problem used in [32] to illustrate the problem of spurious pressure oscillation. This problem was also used in [27] and [38] to demonstrate their methods ability to overcome these spurious pressure oscillations. The cantilever bracket problem (shown in Figure 6.5a) is solved on a unit square $[0,1]^2$. No-flow flux boundary conditions are applied along all sides, the deformation is fixed ($\mathbf{u} = 0$) along the left hand-side ($x = 0$), and a downward traction force, $\mathbf{t}_N \cdot \mathbf{n} = -1$, is applied along the top edge ($y = 1$). The right and bottom sides are traction-free. For this numerical experiment, we set $\Delta t = 0.001$, $h = 1/96$, $\delta = 5 \times 10^{-6}$. The material parameters $\lambda$ and $\mu$ are chosen such that Youngs's modulus, $E = 10^5$ and Poisson's ratio $\nu = 0.4$ and $\alpha = 0.93, c_0 = 0, k = 1 \times 10^{-7}$, values shown in [32] to typically cause locking. The proposed stabilized finite element method yields a smooth pressure solution without any oscillations as is shown in Figure 6.5b.

**6.4. 3D unconfined compression stress relaxation.** In this test, a cylindrical specimen of porous tissue is exposed to a prescribed displacement in the axial direction while left free to expand radially. (Note that the two plates are not explicitly modelled in the simulation, but are realised through displacement boundary conditions.) After loading the tissue, the displacement is held constant while the tissue relaxes in the radial direction due to interstitial fluid flow through the radial boundary. For the special case of a cylindrical geometry [1] found a closed-form analytical solution for the radial displacement $u$ given by

$$\frac{u}{a}(a,t) = \epsilon_0 \left[\nu + (1-2\nu)(1-\nu) \sum_{n=1}^{\infty} \frac{\exp\left(-\alpha_n^2 \frac{Mkt}{a^2}\right)}{\alpha_n^2(1-\nu)^2 - (1-\nu)}\right], \qquad (6.2)$$

where $\alpha_n$ are the solutions to the characteristic equation $J_1(x) - (1-\nu)xJ_0(x)/(1-2\nu) = 0$, where $J_0$ and $J_1$ are Bessel functions, $\epsilon_0$ is the amplitude of the applied axial strain, $a$ is the radius of the cylinder, and $t_g$ is the characteristic time of diffusion (relaxation) $t_g = a^2/Mk$, where $M = \lambda + 2\mu$ is the P-wave modulus of the elastic solid skeleton, and $k$ is the permeability.

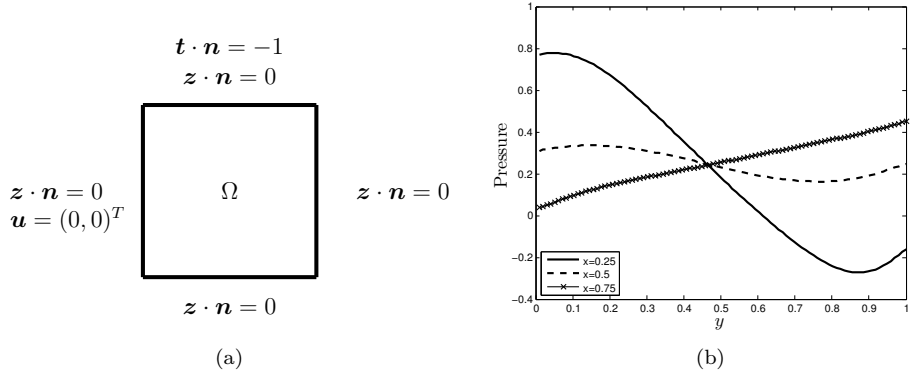

Fig. 6.5: (a) Boundary conditions for the cantilever bracket problem. (b) Pressure solution for the cantilever bracket problem at $t = 0.005$.

The analytical solution available for this test problem describes the displacement of the outer radius which is directly dependent on the amount of mass in the system since the porous medium is assumed to be incompressible and fully saturated. It is therefore an ideal test problem for analyzing the effect that the added stabilization term has on the conservation of mass. In Figure 6.7 we can see that for large values of $\delta$ the numerical solution loses mass faster and comes to a steady state that has less mass than the analytical solution. This is a clear limitation of the method and the stability parameter therefore needs to be chosen carefully. However, for 3D problems $\delta$ can be chosen to be very small so this effect is negligible, as can be seen in Figure 6.7 for a stable value of $\delta = 0.001$.

**7. Conclusion.** The local pressure jump stabilization method [8] is commonly used to solve the Stokes or Darcy equations using piecewise linear approximations for the velocities, and piecewise constant approximations for the pressure variable. The main contribution of this paper has been to extend these ideas to three-field poroelasticity. We have presented a stability result for the discretized equations that guarantees the existence of a unique solution at each time step, and derived an energy estimate which can be used to prove weak convergence of the solution of the discretized system to the solution to the continuous problem as the mesh parameters tend to zero. We also derived an optimal error estimate which includes an error for the fluid flux in its natural $Hdiv$ norm. We have also presented numerical experiments in 2D and 3D that illustrate the convergence of the method, the effectiveness of the method in overcoming spurious pressure oscillations, and the added mass effect of the stabilization term.

REFERENCES

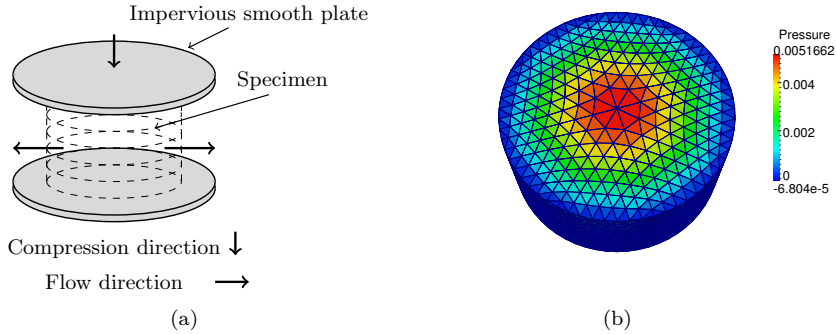

(a) (b)

Fig. 6.6: (a) Sketch of the test problem. The porous medium is being compressed between two smooth impervious plates. The frictionless plates permit the porous medium to expand in order to conserve volume and then to gradually relax as the fluid seeps out radially. (b) Pressure field solution at $t = 5s$, using a mesh with 28160 tetrahedra.

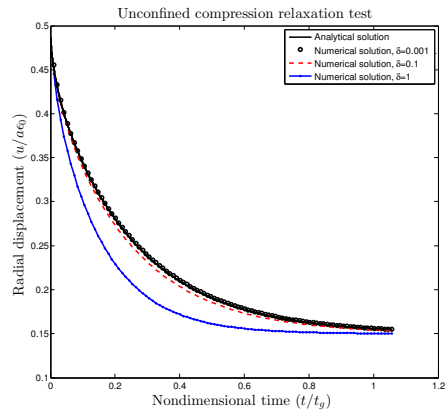

Fig. 6.7: Normalized radial displacement versus normalized time calculated using the analytical solution, and using the proposed numerical method with different values of $\delta$. At $t = 0s$ the radial expansion is half of the axial compression indicating the instantaneous incompressibility of the poroelastic tissue. The final amount of tissue recoil depends on the intrinsic Poisson ratio of the tissue skeleton.

[1] C. G. Armstrong, W. M. Lai, and V. C. Mow, *An analysis of the unconfined compression of articular cartilage*, Journal of Biomechanical Engineering, 106 (1984), pp. 73–165.
[2] I. Babuška, *Error-bounds for finite element method*, Numerische Mathematik, 16 (1971), pp. 322–333.
[3] S. Badia, A. Quaini, and A. Quarteroni, *Coupling Biot and Navier-Stokes equations for modelling fluid-poroelastic media interaction*, Journal of Computational Physics, 228 (2009), pp. 7986–8014.
[4] H. Barucq, M. Madaune-Tort, and P. Saint-Macary, *Some existence-uniqueness results for a class of one-dimensional nonlinear Biot models*, Nonlinear Analysis: Theory, Meth-

ods & Applications, 61 (2005), pp. 591–612.

[5] P. B. BOCHEV AND C. R. DOHRMANN, *A computational study of stabilized, low-order C0 finite element approximations of Darcy equations*, Computational Mechanics, 38 (2006), pp. 323–333.

[6] R. BOER, *Trends in continuum mechanics of porous media*, vol. 18, Springer, 2005.

[7] S. C. BRENNER AND LARKIN R. SCOTT, *The mathematical theory of finite element methods*, vol. 15, Springer, 2008.

[8] E. BURMAN AND P. HANSBO, *A unified stabilized method for Stokes' and Darcy's equations*, Journal of Computational and Applied Mathematics, 198 (2007), pp. 35–51.

[9] V. CAREY, D. ESTEP, AND S. TAVENER, *A posteriori analysis and adaptive error control for operator decomposition solution of coupled semilinear elliptic systems*, International Journal for Numerical Methods in Engineering, 94 (2013), pp. 826–849.

[10] D. CHAPELLE, J.F. GERBEAU, J. SAINTE-MARIE, AND I.E VIGNON-CLEMENTEL, *A poroelastic model valid in large strains with applications to perfusion in cardiac modeling*, Computational Mechanics, 46 (2010), pp. 91–101.

[11] P. G. CIARLET, *The Finite Element Method for Elliptic Problems*, SIAM, classics in applied mathematics ed., 2002.

[12] A. N. COOKSON, J. LEE, C. MICHLER, R. CHABINIOK, E. HYDE, D. A. NORDSLETTEN, M. SINCLAIR, M. SIEBES, AND N. P. SMITH, *A novel porous mechanical framework for modelling the interaction between coronary perfusion and myocardial mechanics*, Journal of Biomechanics, 45 (2012), pp. 850 – 855.

[13] O. COUSSY, *Poromechanics*, John Wiley & Sons Inc, 2004.

[14] D. A. DI PIETRO AND A. ERN, *Mathematical aspects of discontinuous Galerkin methods*, vol. 69, Springer, 2011.

[15] H. C. ELMAN, D. J. SILVESTER, AND A. J. WATHEN, *Finite elements and fast iterative solvers: with applications in incompressible fluid dynamics*, Oxford University Press, USA, 2005.

[16] X. FENG AND Y. HE, *Fully discrete finite element approximations of a polymer gel model*, SIAM Journal on Numerical Analysis, 48 (2010), pp. 2186–2217.

[17] F. GALBUSERA, H. SCHMIDT, J. NOAILLY, A. MALANDRINO, D. LACROIX, H.J. WILKE, AND A. SHIRAZI-ADL, *Comparison of four methods to simulate swelling in poroelastic finite element models of intervertebral discs*, Journal of the Mechanical Behavior of Biomedical Materials, 4 (2011), pp. 1234 – 1241.

[18] P. A. GALIE, R. L. SPILKER, AND J. P. STEGEMANN, *A linear, biphasic model incorporating a Brinkman term to describe the mechanics of cell-seeded collagen hydrogels*, Annals of Biomedical Engineering, 39 (2011), pp. 2767–2779.

[19] J. B HAGA, H OSNES, AND H. P LANGTANGEN, *On the causes of pressure oscillations in low-permeable and low-compressible porous media*, International Journal for Numerical and Analytical Methods in Geomechanics, 36 (2012), pp. 1507–1522.

[20] M. H. HOLMES AND V. C. MOW, *The nonlinear characteristics of soft gels and hydrated connective tissues in ultrafiltration*, Journal of Biomechanics, 23 (1990), pp. 1145–1156.

[21] A-R. A. KHALED AND K. VAFAI, *The role of porous media in modeling flow and heat transfer in biological tissues*, International Journal of Heat and Mass Transfer, 46 (2003), pp. 4989–5003.

[22] J. KIM, H. A. TCHELEPI, AND R. JUANES, *Stability and convergence of sequential methods for coupled flow and geomechanics: Fixed-stress and fixed-strain splits*, Computer Methods in Applied Mechanics and Engineering, 200 (2011), pp. 1591–1606.

[23] P. KOWALCZYK, *Mechanical model of lung parenchyma as a two-phase porous medium*, Transport in Porous Media, 11 (1993), pp. 281–295.

[24] H. LI AND Y. LI, *A discontinuous Galerkin finite element method for swelling model of polymer gels*, Journal of Mathematical Analysis and Applications, 398 (2012), pp. 11–25.

[25] X. G. LI, H. HOLST, J. HO, AND S. KLEIVEN, *Three dimensional poroelastic simulation of brain edema: Initial studies on intracranial pressure*, in World Congress on Medical Physics and Biomedical Engineering, Springer, 2010, pp. 1478–1481.

[26] K. LIPNIKOV, *Numerical methods for the Biot model in poroelasticity*, PhD thesis, University of Houston, 2002.

[27] R. LIU, *Discontinuous Galerkin finite element solution for poromechanics*, PhD thesis, The University of Texas at Austin, 2004.

[28] V. C. MOW, S. C. KUEI, W. M. LAI, AND C. G. ARMSTRONG, *Biphasic creep and stress relaxation of articular cartilage in compression: Theory and experiments*, Journal of Biomechanical Engineering, 102 (1980), pp. 73–84.

[29] M. A. MURAD AND A. F. D. LOULA, *On stability and convergence of finite element approximations of Biot's consolidation problem*, International journal for numerical methods in

engineering, 37 (1994), pp. 645–667.
- [30] P. J. PHILLIPS AND M. F. WHEELER, *A coupling of mixed and continuous Galerkin finite element methods for poroelasticity i: the continuous in time case*, Computational Geosciences, 11 (2007), pp. 131–144.
- [31] ———, *A coupling of mixed and continuous Galerkin finite element methods for poroelasticity ii: the discrete-in-time case*, Computational Geosciences, 11 (2007), pp. 145–158.
- [32] ———, *Overcoming the problem of locking in linear elasticity and poroelasticity: an heuristic approach*, Computational Geosciences, 13 (2009), pp. 5–12.
- [33] R. E. SHOWALTER, *Diffusion in poro-elastic media*, Journal of mathematical analysis and applications, 251 (2000), pp. 310–340.
- [34] R. VERFÜRTH, *A posteriori error estimators for convection-diffusion equations*, Numerische Mathematik, 80 (1998), pp. 641–663.
- [35] M. F. WHEELER AND X. GAI, *Iteratively coupled mixed and Galerkin finite element methods for poro-elasticity*, Numerical Methods for Partial Differential Equations, 23 (2007), pp. 785–797.
- [36] J. A. WHITE AND R. I. BORJA, *Stabilized low-order finite elements for coupled solid-deformation/fluid-diffusion and their application to fault zone transients*, Computer Methods in Applied Mechanics and Engineering, 197 (2008), pp. 4353–4366.
- [37] B. WIRTH AND I. SOBEY, *An axisymmetric and fully 3D poroelastic model for the evolution of hydrocephalus*, Mathematical Medicine and Biology, 23 (2006), pp. 363–388.
- [38] S-Y. YI, *A coupling of nonconforming and mixed finite element methods for Biot's consolidation model*, Numerical Methods for Partial Differential Equations, (2013).
- [39] A. ŽENÍŠEK, *The existence and uniqueness theorem in Biot's consolidation theory*, Aplikace matematiky, 29 (1984), pp. 194–211.